\documentclass[12pt]{amsart}

\usepackage{color}

\usepackage{amsmath}
\usepackage{amsthm}
\usepackage{amssymb}
\usepackage{amscd}
\usepackage{amsfonts}
\usepackage{amsbsy}
\usepackage[noadjust]{cite} %The noadjust option makes (\cite{ref}) space correctly on the left.
\usepackage{epsfig,afterpage}
\usepackage{paralist}
\usepackage{psfrag}
\usepackage{mathrsfs}       %Gives the script font enabled by the invoked command \mathscr
\usepackage{verbatim}

\newtheorem {theorem} {Theorem}%[section]

\newtheorem {proposition} [theorem]{Proposition}
\newtheorem {corollary} [theorem]{Corollary}
\newtheorem {lemma}  [theorem]{Lemma}
\newtheorem {remark} [theorem]{\sc Remark}

\title[Stability of analytic monodromic singularities]{Stability of analytic monodromic singularities without curves of zero angular speed}

\thanks{The author is partially supported by the Agencia Estatal de Investigaci\'on grant number PID2025-172232NB-I00.}

\author[Isaac A. Garc\'ia and Jaume Gin\'e]{Isaac A. Garc\'ia and Jaume Gin\'e}

\address{Departament de Matem\`atica, Universitat de Lleida,
Avda. Jaume II, 69, 25001 Lleida, Spain}

\email{isaac.garcia@udl.cat}
\email{jaume.gine@udl.cat}



\begin{document}

\begin{abstract}
The aim of this work is to give conditions on the parameters of a family of analytic planar vector fields with a fixed Newton diagram and a monodromic singularity in order to guarantee that the coefficient (generalized first Lyapunov constant) of the linear part of the Poincar\'e map has an explicit formula. We revisit using this formalism some monodromic examples appearing in the literature and also we apply the results to other families were, as far as we know, the stability problem was open.
\end{abstract}

\maketitle

\section{Introduction}

In this work, we study families of real-analytic planar differential systems
\begin{equation}\label{DCC-1}
\dot{x}=P(x,y;\lambda), \qquad \dot{y}=Q(x,y;\lambda),
\end{equation}
where the dot denotes differentiation with respect to time. We assume that the vector field depends analytically on the parameters $\lambda\in\mathbb{R}^p$. The family \eqref{DCC-1} defines the vector field $\mathcal{X}=P(x,y;\lambda)\,\partial_x+Q(x,y;\lambda)\,\partial_y$. Throughout the paper, we focus on parameter values for which the origin is a monodromic singularity: there are no characteristic orbits, that is, no orbits approaching the origin with a well-defined tangent direction, and nearby orbits wind around it. We denote by $\Lambda\subset\mathbb{R}^p$ the set of parameters for which the origin is monodromic. Equivalently, for such parameter values, a Poincar\'e first-return map $\Pi$ is defined on a suitable transversal ending at the origin.

The monodromic set $\Lambda$ is commonly characterized using the blow-up process, or $\sigma$-process, developed by Dumortier \cite{Dum}; see also Arnold \cite{Ar}. In \cite{AGR, AGR2}, Algaba and co-authors provide an algorithmic procedure for determining the parameter conditions that characterize $\Lambda$, based on the structure of the Newton diagram $\mathbf{N}(\mathcal{X})$. Related approaches can be found in \cite{GGGrau, BM, Me-monod}.

Determining the stability of a monodromic singularity when the associated linear part of the vector field vanishes is a challenging problem; see, for example, \cite{AGG, GGL, G1, GM}. We also mention \cite{GaGi3}, where the stability problem is studied through the complex separatrices of the singularity.

Desingularization by blow-ups transforms the singularity at the origin into a monodromic polycycle $\Gamma$. The least degenerate case occurs when $\Gamma$ is a periodic orbit, in which case the Poincar\'e map is analytic at the origin. We refer to the monodromic class as ${\rm Mo}^{(p,q)}$ when the periodic orbit $\Gamma=\{\rho=0\}$ is obtained through a $(p,q)$-weighted polar blow-up $(x, y) = (\rho^{p} \cos\varphi, \rho^{q} \sin\varphi)$, see \eqref{wpbu} for details. This class includes among others all monodromic singularities for which the vector field has a nonvanishing associated linear part. ${\rm Mo}^{(p,q)}$ is the simplest monodromic class since the classical Poincar\'e--Lyapunov approach for computing the leading terms in the Taylor expansion of $\Pi$ applies with only minor modifications. For a modern treatment, see \cite{GaGi3,GaGi5} and the book \cite{RS}. Our interest lies instead in the more difficult cases, in which $\Gamma$ becomes a monodromic polycycle, hence contains singularities of the vector field. Recent work in this direction includes \cite{Ga-Gi, GaGi2}.

It is known, see Il'yashenko's work \cite{Ilyashenko}, that the Poincar\'e map $\Pi$ is a semiregular map that has a linear part and can be written in the form
\[
\Pi(x;\lambda)=\eta_1(\lambda)x+o(x),
\]
where the leading coefficient $\eta_1(\lambda)>0$ is called the first generalized Poincar\'e--Lyapunov quantity, see \cite{Ga-Ma-Ma}. This key result was stated without proof in \cite{AA} and was later proved in \cite{Me-2}, see also \cite{GaGi4} for a different proof under some technical restrictions.

If $\eta_1(\lambda)\neq 1$, then the origin is a structurally stable focus, and its stability is determined by the value of $\eta_1(\lambda)$. When $\eta_1(\lambda)=1$, additional terms in the asymptotic Dulac expansion of $\Pi$ are required to determine stability, see for instance \cite{Me-2,Me-3}. We do not address this case in the present work.

Outside the class ${\rm Mo}^{(p,q)}$, computing $\eta_1$ is cumbersome and is generally possible only under certain nondegeneracy assumptions. Such computations rely on a method proposed by H. Dulac in his classical monograph \cite{Dulac}; see also \cite{Ilyashenko}. In this work, we contribute to the stability analysis by deriving an explicit formula for $\eta_1$ for monodromic singularities without curves of zero angular speed, see Theorem \ref{cor-calculo-eta1}.
\newline

The work is structured as follows. The main results are stated in section \ref{S-main} and culminates in Theorem \ref{cor-calculo-eta1} where we describe a method to determine the stability of a monodromic singularity without curves of zero angular speed under some conditions previously stated in Theorem \ref{teo-Fund-G-jet}. In section \ref{gm-cw} we prove some results on the existence of analytic solutions of singular homogebeous linear differential equations needed in the applications of our results. Section \ref{sec-V} is devoted to explain the relation of our framework with the existence of an inverse integrating factor. In the section  \ref{sec-ex} we give some nontrivial illustrative examples. We relegate to the last section the proofs of the main results.

\section{Main results} \label{S-main}

Before the statement of the main results we are going to introduce some needed background. Hereinafter we omit the dependency on the parameters $\lambda$ in the expressions except when it is relevant for discussing.

We denote by $\mathbf{N}(\mathcal{X})$ the Newton diagram of $\mathcal{X}$, consisting of edges that connect the positive coordinate axes. For a detailed account of this construction, see \cite{Bruno}; see also \cite{AGR}. The endpoints of each edge of $\mathbf{N}(\mathcal{X})$ lie in $\mathbb{N}^2$, and its slope is $-p/q$, where $p,q\in\mathbb{N}$ are coprime. We denote by $W(\mathbf{N}(\mathcal{X}))\subset\mathbb{N}^2$ the set of all such weight pairs $(p,q)$. For each weight pair, we associated the weighted polar blow-up
\begin{equation}\label{wpbu}
x=\rho^p\cos\varphi,\qquad y=\rho^q\sin\varphi.
\end{equation}
In the coordinates $(\rho, \varphi)$, and after removing a power of $\rho$ that appears as common factor, the vector field $\mathcal{X}$ becomes the polar vector field $\mathcal{Z} = \Theta(\varphi, \rho) \partial_\varphi + R(\varphi, \rho) \partial_\rho$ defined on the cylinder $C =  \{ (\varphi, \rho) \in \mathbb{S}^1 \times \mathbb{R} :  0 \leq |\rho| \ll 1 \}$ where $\mathbb{S}^1 = \mathbb{R}/ (2 \pi \mathbb{Z})$. By monodromy, $R(\varphi, 0)=0$ and $\Theta(\varphi, 0) \geq 0$ holds. The set of {\it characteristic directions} at the origin of $\mathcal{X}$ is defined as
$\Omega_{pq} = \{ \varphi^* \in \mathbb{S}^1 : \Theta(\varphi^*, 0) = 0 \}$. The set $\{\rho=0 \}$ is invariant under the flow of $\mathcal{Z}$ and is either a periodic orbit or a polycycle, according to whether $\Omega_{pq} = \emptyset$ or $\Omega_{pq} \neq \emptyset$, respectively. In particular, $\Pi$ is analytic when $\Omega_{pq} = \emptyset$ that corresponds to the monodromic class ${\rm Mo}^{(p,q)}$. We call
$\Theta^{-1}(0) = \left\{(\varphi,\rho)\in C: \Theta(\varphi,\rho)=0,\ \rho\geq 0\right\}$ the set of zero angular speed. We also define the $(p,q)$-critical parameters as the elements of the subset $\Lambda_{pq}\subset\Lambda$ of the monodromic parameter space for which the vector field has local curves of zero angular speed:
\[
\Lambda_{pq}
=
\left\{
\lambda\in\Lambda:
\Theta^{-1}(0)\setminus\{\rho=0\}\neq\emptyset
\right\}.
\]
We emphasize that $\Lambda_{pq}=\emptyset$ whenever $\Omega_{pq}=\emptyset$, although the converse need not hold. Notice that when the parameters lie in $\Lambda \setminus \Lambda_{pq}$, the set $\Theta^{-1}(0)$ is discrete and only contains the singularities on $\rho=0$ of $\mathcal{Z}$. It is also worth to emphasize that, in \cite{Ga-Li-Ma-Ma}, it is proved that $\Lambda_{11} = \emptyset$ when $\mathcal{X}$ is polynomial vector field defined by the sum of two homogeneous vector fields. We refer the reader to the appendix of \cite{GaGi4} for a detailed explanation about the explicit computation of the set $\Lambda_{pq}$ of critical parameters associated to  curves of zero angular speed.

\begin{remark}
{\rm We want to emphasize that the existence or not of curves of zero angular speed may depend on the chosen weights $(p,q) \in W(\mathbf{N}(\mathcal{X}))$. We show this behaviour with family \eqref{ejemplo3-JDE-Armengol} with monodromic parameter set $\Lambda = \{ (a, b) \in \mathbb{R}^2 : b > 1/4 \}$ analyzed in section \ref{sec-ex}. This family has $W(\mathbf{N}(\mathcal{X})) = \{ (1,1), (1,2) \}$. First we see that $\Lambda_{11} = \emptyset$ because, using polar coordinates $(\varphi, \rho)$, the equation $\Theta(\varphi, \rho)$ is polynomial of degree 2 in $\rho$ with discriminant negative for any $\varphi$ in a neighborhood of $\Omega_{11} = \{0, \pi\}$ and $b \in \Lambda$. On the other hand, if now $(\varphi, \rho)$ denote (1,2)- weighted polar coordinates, then $\Theta(\varphi, \rho)$ is again a quadratic polynomial in $\rho$ but its discriminant is positive when $\varphi$ lies in a neighborhood of $\Omega_{12} = \{\pi/2, 3 \pi/2\}$ and $b \in \Lambda$ provided $|a| > 4 b$. In summary, we have $\Lambda_{12} \neq \emptyset$. }
\end{remark}

\subsection{The function $\mathcal{G}$ and its fundamental equation}

We consider the rectangle $\mathcal{R} = [0, 2\pi] \times I$, where $I \subset \mathbb{R}$ is a neighborhood of the origin. We seek for solutions $\mathcal{G} : \mathcal{R} \to \mathbb{R}$ of the linear partial differential equation
\begin{equation}\label{edp-G}
\mathcal{Z}(\mathcal{G}) = \mathcal{G} \; {\rm div}(\mathcal{Z}),
\end{equation}
where ${\rm div}(\mathcal{Z}) = \partial_\varphi \Theta + \partial_\rho R$ is the divergence of $\mathcal{Z}$. Notice that \eqref{edp-G} is equivalent to the divergence-free condition ${\rm div}(\mathcal{Z}/\mathcal{G}) \equiv 0$ valid in $\mathcal{R} \setminus \mathcal{G}^{-1}(0)$.

\begin{remark}\label{Rem-G-cero-linear0}
{\rm We will also use the rescaled vector field $\mathcal{Z}^* = \mathcal{Z}/\Theta = \partial_\varphi + \mathcal{F}(\varphi, \rho) \partial_\rho$ with $\mathcal{F}=R/\Theta$ that is analytic in $C \setminus \Theta^{-1}(0)$. We define the function $\mathcal{G}^* : \mathcal{R} \setminus \Theta^{-1}(0) \to \mathbb{R}$ satisfying the condition ${\rm div}\left( \mathcal{Z}^*/ \mathcal{G}^* \right) \equiv 0$ in $\mathcal{R}$ except in $\Theta^{-1}(0) \cup (\mathcal{G}^*)^{-1}(0)$. This the divergence free condition is equivalent, in $\mathcal{R} \setminus \{ \Theta^{-1}(0) \cup (\mathcal{G}^*)^{-1}(0) \}$, to the linear partial differential equation
\begin{equation}\label{edp-G*}
\mathcal{Z}^*(\mathcal{G}^*) = \mathcal{G}^* \, {\rm div}\left(\mathcal{Z}^*\right).
\end{equation} }
\end{remark}

\begin{remark}\label{Rem-G-cero-linear}
{\rm Even in the simplest case $\Theta^{-1}(0) = \emptyset$, equation \eqref{edp-G*} can have solutions $\mathcal{G}^*$ analytic in $\mathcal{R}$ but not in $C$ since they are not $2 \pi$-periodic in $\varphi$. An example is given by the linear vector field $\mathcal{X} = (-y + \lambda x) \partial_x + (x + \lambda y) \partial_y$ with $\lambda \neq 0$. In polar coordinates, $\mathcal{Z} = \mathcal{Z}^* = \partial_\varphi + \lambda \rho \partial_\rho$ and \eqref{edp-G} has the family of solutions  $\mathcal{G}^*(\varphi, \rho) = \rho - \rho_0 \exp(\lambda \varphi)$ with $\rho_0 \in \mathbb{R}$. Nevertheless we see that its zero-set $(\mathcal{G}^*)^{-1}(0)$ is a logarithmic spiral in $C$ when $\rho_0 > 0$ that corresponds to an spiral orbit of the focus approaching or leaving $\rho=0$. Therefore, $(\mathcal{G}^*)^{-1}(0) \cap C \backslash \{\rho=0\} \neq \emptyset$ although clearly $(\mathcal{G}^*)^{-1}(0) \cap \mathcal{R} = \emptyset$ taking the length of $I$ sufficiently small since $0 \leq \varphi \leq 2 \pi$ in $\mathcal{R}$.  This behaviour is also  reproduced under the parameter restriction $\Lambda \setminus \Lambda_{pq}$ as it is stated Lemma \ref{structure-G=0}. }
\end{remark}

\begin{theorem} \label{teo-Fund-GPi2}
We restrict to the parameter space $\Lambda \setminus \Lambda_{pq}$ and, by a rotation, we take $0 \not\in \Omega_{pq}$. Then there are real analytic functions $\mathcal{G}(\varphi, \rho)$ defined on $\mathcal{R} \setminus \Theta^{-1}(0)$ such that ${\rm div}\left( \mathcal{Z}/ \mathcal{G} \right) \equiv 0$. The function $\mathcal{G}$ is unique once the initial condition $\mathcal{G}(0, \rho)$ analytic in $I$ is fixed. Additionally, the fundamental equation
\begin{equation}\label{G-fundamental2}
\mathcal{G}(2 \pi, \Pi(\rho_0)) \, \Theta(0, \rho_0) = \mathcal{G}(0, \rho_0) \, \Theta(0, \Pi(\rho_0)) \, \Pi'(\rho_0).
\end{equation}
is satisfied where $\Pi$ is the Poincar\'e map of the monodromic singularity.
\end{theorem}

\begin{remark}\label{Rem-CVAN}
{\rm In general, Theorem \ref{teo-Fund-GPi2} does not work when the parameters lie in $\Lambda_{pq}$. The reason is the existence of curves $\Theta^{-1}(0) \subset C$ of zero angular speed that prevent the construction of a global $\mathcal{G}$ in $\mathcal{R} \setminus \{\rho=0\}$ by the characteristics method in the partial differential equation \eqref{edp-G}. This global solution is needed in the proof of equation \eqref{G-fundamental2}. }
\end{remark}

\subsection{A method to determine the stability in $\Lambda \setminus \Lambda_{pq}$} \label{sec-Algor}

The goal now is to compute, under some technical conditions, the coefficient $\eta_1$ of the linear part of $\Pi$. First we will need to perform, associated to each $\varphi^* \in \Omega_{pq}$, a suitable blow-up to the partial differential equation \eqref{edp-G} as described in the next result.

\begin{theorem} \label{teo-Fund-G-jet}
In $\Lambda \setminus \Lambda_{pq}$ and with $0 \not\in \Omega_{pq}$, we take the counterclockwise sense flow rotation. Let $\mathcal{G}$ be any function described in Theorem \ref{teo-Fund-GPi2} such that $\mathcal{G}(0, \rho)$ has order $m \geq 2$ at $\rho=0$. Then the restriction of $\mathcal{G}$ to $I_{pq} = [0, 2 \pi] \setminus \Omega_{pq}$ admits the convergent Taylor expansion
\begin{equation}\label{Taylor-G}
\mathcal{G}(\varphi, \rho) \Big|_{I_{pq}} = g_m(\varphi) \, \rho^m + o(\rho^m).
\end{equation}
For each $\varphi^* \in \Omega_{pq}$ there are nonzero constants $C_{\pm}$ such that $g_m(\varphi) = C_{-} \, \hat{g}_m(\varphi)$ when $\varphi \lesssim \varphi^*$ while $g_m(\varphi) = C_{+} \, \hat{g}_m(\varphi)$ if $\varphi \gtrsim \varphi^*$ for the $\mathcal{G}$ with $\mathcal{G}(0, \rho) = \hat{g}_m(0) \rho^m + \cdots$ where $\hat{g}_m(\varphi)$ is a solution of the linear singular differential equation
\begin{equation}\label{EDO-GM}
G_0 \hat{g}'_m + ((m-1) R_1 - G'_0) \hat{g}_m = 0.
\end{equation}
Assume moreover that $\hat{g}_m$ is analytic at $\Omega_{pq}$ and that, for each $\varphi^* \in \Omega_{pq}$, there are $\alpha > 0$ and $r > 0$ such that the blow-up $\varphi-\varphi^* = \rho^\alpha u$  yields
\begin{eqnarray*}
R(\varphi^* + \rho^\alpha u, \rho) &=& \rho^r Q(u) + o(\rho^{r}), \\
\Theta(\varphi^* +\rho^\alpha u, \rho) &=& \rho^{r+\alpha-1} P(u) + o(\rho^{r+\alpha-1}), \\
a(\varphi^* +\rho^\alpha u, \rho) &=& \rho^{\nu} B(u) + o(\rho^{\nu}),
\end{eqnarray*}
for some $\nu \geq 0$ and certain analytic function $a(\varphi-\varphi^*, \rho)$ with $a(\varphi-\varphi^*, 0) = \hat{g}_m(\varphi-\varphi^*)$. We define the function $D(u) = P(u)-\alpha u Q(u)$ and we assume that there exist constants $c_D, c_B > 0$ such that, for every $u\in\mathbb{R}$, the lower bounds
\begin{equation}\label{eq:positiveDB}
 D(u)\ge c_D(1+|u|^{1/\alpha})^{r+\alpha-1}, \qquad B(u)\ge c_B(1+|u|^{1/\alpha})^\nu
\end{equation}
hold. Then
\begin{equation}\label{eq:Cconnection}
 \frac{C_-}{C_+}=\exp\left(\int_{-\infty}^{+\infty}\Xi(u)\,d u\right),
\end{equation}
where the integral converges absolutely and
\begin{equation}\label{eq:Xi1}
\Xi(u)=\frac{D(u)B'(u)+B(u)\bigl((m+\nu-r)Q(u)+\alpha uQ'(u)-P'(u)\bigr)}{B(u)D(u)}.
\end{equation}
\end{theorem}

Since \eqref{EDO-GM} is singular at $\Omega_{pq}$, the existence of a $\hat{g}_m$ analytic at $\Omega_{pq}$ is not guarantee although it is necessary in order to apply Theorem \ref{teo-Fund-G-jet}. Indeed, sometimes this analyticity is only obtained for concrete values of $m$. We will write in the next section a detailed analysis about the analyticity of the solutions of \eqref{EDO-GM}.

On the other hand, the bound on $D(u)$ in \eqref{eq:positiveDB} implies in particular $D(u) > 0$. This is the one-chart desingularization condition  since, after the blow-up and rescaling the time removing the factor $\rho^{r-1}$, the polar vector field $\mathcal{Z}$ is written as $\dot{\rho} = Q(u) \rho + O(\rho^2)$ and $\dot{u} = D(u) + O(\rho)$, hence the blow-up removes all the singularities of the reduced vector field on the finite part of the exceptional divisor.

\begin{remark}
{\rm The analyticity of $\hat{g}_m$ at $\Omega_{pq}$ for some $m \geq 2$ depends strongly on the $(p,q)$-weighted polar coordinates $(\rho, \varphi)$ used to compute the linear singular differential equation \eqref{EDO-GM}. See this behaviour in the analysis of the forthcoming family \eqref{Ej2-Nonlinearity}. }
\end{remark}

\begin{remark} \label{Remar-funcion-a}
{\rm If $\hat{g}_m$, with $m \geq 2$, is a solution of the linear singular differential equation \eqref{EDO-GM} analytic at $\varphi^*$ such that $\hat{g}_m(\varphi^*) \ne 0$, one may take the function $a(x,\rho)$ of Theorem \ref{teo-Fund-G-jet} with $x=\varphi-\varphi^*$ to be
$a(x, \rho) = \hat{g}_m(x)$, so that $\nu=0$ and $B$ is a constant. On the contrary, if $\hat{g}_m(\varphi^*) = 0$ then (see the proof of Theorem \ref{teo-Fund-G-jet}) $\hat{g}_m(x) \rho^m$ has zeros along the orbit  and one instead chooses a function $a(x,\rho)$ that is non-zero for $\rho>0$ and has, of course, $a(x,0)=\hat{g}_m(x)$. In this case we will take $a(x, \rho) = \hat{g}_m(x) + c_* \rho^{\kappa \alpha}$ where $c_* > 0$ is a constant and $\kappa \geq 0$ is the vanishing multiplicity of $\hat{g}_m$ at $\varphi^*$. The reason of the choice of the exponent $\kappa \alpha$ is because if we take another non-negative exponent less than it then $B(u) = c_*$ and $\nu>0$ so that the bound on $B(u)$ given in \eqref{eq:positiveDB} is never possible. }
\end{remark}

\begin{remark} \label{Remar-cotas}
{\rm  Regarding the lower bounds \eqref{eq:positiveDB}, we define $\beta(u) = (1+|u|^{a})^{b}$ with constants $a, b > 0$ and we note that if $F(u) >0$ and continuous in $\mathbb{R}$ and $F(u)/\beta(u)$ has finite limit when $|u| \to \infty$ then $\inf_{\mathbb{R}} F/\beta > 0$ and consequently there is $c > 0$ such that $F(u) \geq c \beta(u)$.  }
\end{remark}

\begin{remark} \label{Remar-gm-Cw}
{\rm We think that the condition on the analyticity of $\hat{g}_m$ at $\Omega_{pq}$ stated in Theorem \ref{teo-Fund-G-jet} can be relaxed to $\hat{g}_m$ to be $C^k(\Omega_{pq})$ for some $k \in \mathbb{N} \cup \{ \infty \}$. We will analyze this possibility in a future work.  }
\end{remark}

We describe below an alternative method to compute the ratio $C_+/C_-$. Let $\alpha$ be the exponent of the admissible blow-up $x= \rho^\alpha u$ associated to the characteristic direction $\varphi^*$ according to Theorem \ref{teo-Fund-G-jet}. Let $\kappa \geq 0$ be the multiplicity of the analytic function $\hat{g}_m$ at $\varphi^*$.

We have $\mathcal{G}(\varphi^* + x, \rho) = g_m(\varphi^*+x) \rho^m + o(\rho^m)$ where $g_m(\varphi^*+x) = C_{+} K x^\kappa + o(x^\kappa)$ or $g_m(\varphi^*+x) = C_{-} K x^\kappa + o(x^\kappa)$ according to whether $x \gtrsim 0$ or $x \lesssim 0$, respectively. Then $\mathcal{G}(\varphi^* + \rho^\alpha u, \rho) = g_m(\varphi^*+ \rho^\alpha u) \rho^m  + \cdots$ where the dots are other terms. We have no control on these terms because we do not know the behaviour of $\hat{g}_j(\varphi^*+x)$ at $x=0$ for $j > m$. Anyway we may guess that these dot terms have no power in $\rho$ less than $\rho^{m + \alpha \kappa}$. Notice that this assumption is equivalent to assume the existence of a function $h(u)$ such that $\mathcal{G}(\varphi^* + \rho^\alpha u, \rho) = h(u) \rho^{m + \alpha \kappa} + o(\rho^{m + \alpha \kappa})$.

\begin{proposition}\label{prop-perfil}
Under the conditions of Theorem \ref{teo-Fund-G-jet}, for each $\varphi^* \in \Omega_{pq}$, let $\kappa \geq 0$ be the multiplicity of $\hat{g}_m$ at $\varphi^*$. We define $h(u)$ to be a nontrivial solution of the linear homogeneous differential equation
\begin{equation}\label{profile-edo}
D(u) h'(u) = \big( P'(u) + (r-m-\alpha \kappa) Q(u)- \alpha u Q'(u) \big) h(u),
\end{equation}
where $D(u) = P(u)-\alpha u Q(u)$. Assume moreover that $L := \lim_{u\to \infty}h(u)/h(-u)$ exists, is finite, and is nonzero. If the compatibility condition
\[
\int_{-\infty}^{+\infty}\left[\frac{B'}{B}+(\nu-\alpha\kappa)\frac{Q}{D}\right] du = 0
\]
holds then $|L|=C_+/C_-$.
\end{proposition}

Sufficient conditions for the existence of $L$ are given in Proposition \ref{prop-perfil-2}.

In the next result we show how, when Theorem \ref{teo-Fund-G-jet} can be applied to each $\varphi^* \in \Omega_{pq}$, we solve the stability problem at first order for the origin of a monodromic singularity, namely we compute the linear term of the Poincar\'e map.

\begin{theorem} \label{cor-calculo-eta1}
In $\Lambda \setminus \Lambda_{pq}$ with $0 \not\in \Omega_{pq}$, assume that Theorem \ref{teo-Fund-G-jet} can be applied to each $\varphi_j^* \in \Omega_{pq}$ and we set
$$
\Xi_j(u) = \frac{D_j(u)B'_j(u)+ B_j(u) \bigl((m + \nu_j- r_j) Q_j(u) + \alpha_j u Q_j'(u)- P_j'(u) \bigr)}{B_j(u) D_j(u)}.
$$
Then the Poincar\'e map is $\Pi(\rho) = \eta_1 \rho + o(\rho)$ with
\begin{equation}\label{eta1}
\eta_1 = \left( C_m \frac{\hat{g}_m(0)}{\hat{g}_m(2 \pi)} \right)^{\frac{1}{m-1}} = C_m^{1/(m-1)} \, \ \exp\left( PV \int_{0}^{2 \pi} \frac{R_1(\varphi)}{G_0(\varphi)} d \varphi   \right)
\end{equation}
where $C_m = \prod_j \exp\left( \int_{\mathbb{R}} \Xi_j(u) \, du \right)$.
\end{theorem}

\begin{remark}
{\rm The last equality in \eqref{eta1} holds because, by the hypothesis of Theorem \ref{teo-Fund-G-jet}, $\hat{g}_m$ is analytic at $\Omega_{pq}$ and then $R_1/G_0$ has at most a simple pole at $\Omega_{pq}$, hence the Cauchy principal value that appears in \eqref{eta1} exists, see for more details Corollary \ref{teo-eta1-m}. An example where this principal value does not exist is given in the forthcoming family \eqref{Ej2-Nonlinearity}. }
\end{remark}

\section{On the analyticity of $\hat{g}_m$} \label{gm-cw}

The leading coefficient $\hat{g}_m$ of $\mathcal{G}$ is a nontrivial solution of the homogeneous linear differential equation \eqref{EDO-GM}. Only analytic solutions with $m > 1$ of \eqref{EDO-GM} are guarantee in each connected open component $I_{pq}$ since \eqref{EDO-GM} is singular at $\Omega_{pq}$. In each connected component $I_k$ of $I_{pq}$ with $k=1, \ldots, \# \Omega_{pq} + 1$, we have up to a multiplicative constant (which may depend on the parameters $\lambda$ of the family) that
\begin{equation}\label{vm-explicit}
\hat{g}_m(\varphi) = G_0(\varphi) \exp\left( (1-m) \, \mathcal{P}(\varphi) \right)
\end{equation}
with $\mathcal{P}' = R_1/G_0$.

Let us analyze a very simple linear differential equation singular at $x=0$. We consider $x y' - y = 0$ and we separate variables $y'/y = 1/x$ and after integrating $\log |y| = \log |x| + c$ we get its general solution
$$
y(x) = C |x| = \left\{ \begin{array}{lll} C x & \mbox{if} & x \geq 0, \\ -C x & \mbox{if} & x < 0, \end{array} \right.
$$
with arbitrary constant $C$ obtained via this process. We see that the only analytic solution at $x=0$ inside this family is just the trivial one $y(x) = 0$. Nevertheless there are other analytic solutions at $x=0$ given by $y(x) = C x$ for any $C$ that are special in the sense that they are not contained in the initial integrating process. Indeed, these special analytic solutions corresponds to choice of constants in the piecewise solutions
$$
y(x) = \left\{ \begin{array}{lll} C_{1} x & \mbox{if} & x \geq 0, \\ - C_{2} x & \mbox{if} & x < 0, \end{array} \right.
$$
given $C := C_{1} =- C_{2}$.

We will use the following general lemma on analytic homogeneous singular differential equations. Essentially it indicates that the following (i), (ii) and (iii) are equivalent, for homogeneous singular analytic  linear differential equations $A(x) y' + B(x) y = 0$:  (i) there is a formal solution; (ii) there is an  analytic solution; (iii) $-B/A$ has (at most) one simple pole and its residue is a non-negative integer $\mathbb{N}_0$. In the following result we refer to analytic special solution an analytic solution that is not contained in the separation of variables integrating process.

\begin{lemma} \label{lemma-edo-homogenea-formal}
Let $A(x) y' + B(x) y = 0$ be a singular differential equation whose coefficients are real analytic in a neighborhood of $x=x^*$ and $A(x^*)=0$. Then any formal solution at $x=x^*$ converges. Moreover, there is a nontrivial real analytic solution $y(x) \not\equiv 0$ if and only if $x=x^*$ is a regular singularity (that is, $B/A$ has at most has a simple pole at $x=x^*$) and $r = {\rm Res}(-B/A, x^*) \in \mathbb{N}_0$ is just the order at $x=x^*$ of that solution. If moreover $y(x)$ is not an analytic special solution then $r \in 2 \mathbb{N}_0 = \{0, 2, 4, \ldots\}$.
\end{lemma}

\begin{remark}
{\rm The analysis of the formal and analytic solutions in the analytic singular non-homogeneous linear equations $A(x) y' + B(x) y = C(x)$ with $A(0)=0$ and $C \not\equiv 0$ is cumbersome. There are formal solutions that does not converge like $y(x) = \sum_{n \geq 0} n! x^{n+1}$ for $x^2 y' - y = -x$. Another behaviour that can happen is that the equation has analytic non-trivial solution but the associated homogeneous one has not like in $x^2 y' - y = -1$ with $y(x) = -1$.   }
\end{remark}

The analyticity and the order of a solution $\hat{g}_m$ at $\varphi^*$ is controlled applying Lemma \ref{lemma-edo-homogenea-formal} to the particular case of the singular differential equation \eqref{EDO-GM} that defines $\hat{g}_m$ with the identification $A = G_0$, $B = (m-1) R_1 - G'_0$, and $x^* = \varphi^*$. Hence we can write the next result.

\begin{corollary} \label{teo-eta1-m}
The nontrivial solution $\hat{g}_m$ of the singular differential equation \eqref{EDO-GM} for some $m \geq 2$ is analytic in $[0, 2\pi]$ if and only if the function $F := ((m-1) R_1 - G'_0)/G_0$ has, for each $\varphi^* \in (0, \pi)$, at most a simple pole and moreover ${\rm Res}(-F,\varphi^*) \in 2 \mathbb{N}_0$.
\end{corollary}

We notice that the multiplicities $2 k$ and $\ell$ of $G_0$ and $R_1$ at any $\varphi^* \in \Omega_{pq}$ generically satisfy the condition $\ell \geq 2k-1$ and is therefore related to the applications of Corollary \ref{teo-eta1-m}. In particular the $F$ of Corollary \ref{teo-eta1-m} has a simple pole at  $\varphi^*$ when $\ell = 2k-1$ and below we see that this is the generic case.

\begin{lemma}\label{lemma-multipl}
In the monodromic scenario, if $G_0$ vanishes at $\varphi^* \in \Omega_{pq}$ with multiplicity $2 k$ then $R_1$, generically, vanishes at $\varphi^*$ with  multiplicity at least $2k-1$.
\end{lemma}

Examples of non-generic monodromic points in the sense of Lemma \ref{lemma-multipl} are, for example, those in case (b) of Theorem 13 in \cite{AGR2}.

\begin{remark}
{\rm The monodromic class $S_{k \omega}$ was introduced in the work \cite{Ga-Ma-Ma}. It is defined as the set of planar vector fields $\mathcal{X}$ analytic at the origin with first non-zero jet of order $k$, the orbits of $\mathcal{X}$ rotates in the positive sense (counterclockwise), and each singular point on the first polar blow-up polycycle desingularizes into two well-posed hyperbolic saddles. Using polar coordinates the class $S_{k \omega}$ is characterized by $0 \not\equiv G_0(\varphi) \geq 0$ and for any $\varphi^* \in \Omega_{11} = G_0^{-1}(0)$ the following conditions hold: (a) $R_1(\varphi^*) = 0$; (b) $G_{1}(\varphi^*) = 0$; (c) $[(G''_{0}-2 R'_{1}) G''_{0}](\varphi^*) > 0$ and (d) $[(G'_{1}- R_{2})^2 -2 (G''_{0} - 2 R'_{1}) G_{2}](\varphi^*) < 0$.

In \cite{Ga-Ma-Ma} it is proved that the origin of $\mathcal{X} \in S_{k \omega}$ has a monodromic singularity at the origin and that $\eta_1$ is given by \eqref{eta1} with $C_m=1$. We observe that sometimes $(1,1) \not\in W(\mathbf{N}(\mathcal{X}))$ but the expression for $\eta_1$ is still valid, see \cite{GaGi5}. We claim that, since $G_{0}(\varphi^*) = G_{0}'(\varphi^*) = 0$, under the conditions (a) and (b) it follows that the function $F$ has at most a simple pole at $\varphi^*$ and moreover ${\rm Res}(-F,\varphi^*) = 2(G_{0}''(\varphi^*) - (m-1) R_{1}'(\varphi^*)) / G_{0}''(\varphi^*) > 0$ for $m=3$ when (c) holds. }
\end{remark}

\section{Inverse integrating factors} \label{sec-V}

Let $V^* : C \to \mathbb{R}$ be an inverse integrating factor of $\mathcal{Z}^* = \partial_\varphi + \mathcal{F}(\varphi, \rho) \partial_\rho$, that is,  ${\rm div}\left( \mathcal{Z}^*/ V^* \right) \equiv 0$. In particular, $V^*$ is $2 \pi$-periodic in the variable $\varphi$. The inverse integrating factor has played a significant role in the center problem, see \cite{EncPer,GP,GGR,Ga-Ma}.
Adapting Theorem 3 of \cite{GaGiGr} to our setting it follows easily that
\begin{equation}\label{fundamental}
V^*(0, \Pi(\rho_0)) = V^*(0, \rho_0) \, \Pi'(\rho_0)
\end{equation}
holds when $\Omega_{pq} = \emptyset$, that is, when $\{\rho=0\}$ is a periodic orbit of $\mathcal{Z}^*$. Later, in Theorem 3 of \cite{GaGi4} the ideas of \cite{GaGiGr} where improved and it was  proved that indeed \eqref{fundamental} holds in the restricted parameter space $\Lambda \backslash \Lambda_{pq}$, hence allowing that $\{\rho=0\}$ be a polycycle but keeping $\Theta^{-1}(0) \backslash \{\rho=0\} = \emptyset$.

\begin{remark}\label{Rem-G-V}
{\rm Despite the condition ${\rm div}\left( \mathcal{Z}^*/ \mathcal{G}^* \right) \equiv 0$ that defines $\mathcal{G}^*$, we emphasize that $\mathcal{G}^*$ is not an inverse integrating factor of $\mathcal{Z}^*$ except when $\mathcal{G}^*$ is defined in $C = \mathbb{S}^1 \times I$ instead of $\mathcal{R}$. In this sense the function $\mathcal{G}^*$ is a generalization of the concept inverse integrating factor $V^*$.}
\end{remark}

Taking the cylinder $C = \mathbb{S}^1 \times I$ with $I = I^- \cup \{0\} \cup I^+$ and $I^+$ is a sufficiently small positive half-neighborhood of the origin, the function $G : I^+ \to \mathbb{R}$ given by
\begin{equation}\label{Function-G}
G(r) = \int_0^{2 \pi} \frac{\mathcal{F}(\varphi, r)}{V^*(\varphi, r)} d \varphi
\end{equation}
is well defined when the parameters are restricted to $\Lambda \backslash \Lambda_{pq}$, see \cite{GaGi4}. The following result is an improvement of Theorem 4 and Theorem 7 of \cite{GaGi4}.

\begin{theorem} \label{lemma-fundamental}
Let $V^*(\varphi, \rho) = v_m(\varphi) \rho^m + \cdots$ be a Laurent expansion at $\rho=0$ of an inverse integrating factor of $\mathcal{Z}^*$ with $m \in \mathbb{Z}$. Then the Poincar\'e map $\Pi$ of the singularity at the origin restricted to the parameter space $\Lambda \backslash \Lambda_{pq}$ admits a formal expansion and has the following structure:
\begin{itemize}
  \item[(i)] If $m \leq 0$ then the origin is a center.

  \item[(ii)] If $m > 1$ then $\Pi(\rho) = \rho + \eta_m \rho^m (1+ O(\rho))$.

  \item[(iii)] If $m=1$ and $\eta_1 = 1$ then the origin is a center.

  \item[(iv)] $G(r) = \mathfrak{g} \in \mathbb{R}$ is a constant in $I^+$. Moreover, taking $v_m(0)=1$ without loss of generality, then $\log\eta_1 = \mathfrak{g}$ when $m=1$ and $\eta_m = \mathfrak{g}$ if $m > 1$.
\end{itemize}
\end{theorem}

\begin{remark}\label{Rem-G-V2}
{\rm The condition $v_m(0)=1$ needed in the statement (iv) of Theorem \ref{lemma-fundamental} is very important to get the correct value of $\eta_m$ using the function $G$ in \eqref{Function-G}, hence to the determine the stability of the monodromic singularity. This condition is missing in Theorem 7 of \cite{GaGi4}. }
\end{remark}

\section{Examples} \label{sec-ex}

In this section we show how our theory works in nontrivial examples. In all the examples, when $0 \in \Omega_{pq}$ but $\pi/2 \not\in \Omega_{pq}$, instead of a rotation we do the linear change of coordinates $(x,y) \mapsto (y, -x)$. The origin of the resulting vector field keeps the stability and the flow sense rotation (since the Jacobian determinant of the change is positive) but now $0 \not\in \Omega_{pq}$.

\subsection{An example from \cite{Ga-Li-Ma-Ma} in the $\mathcal{G}$-monodromic class}

In \cite{Ga-Li-Ma-Ma} it is shown that the origin of family
\begin{equation}\label{Ej2-Nonlinearity}
\dot{x}= y(\alpha x^2 + b x y+ c y^2),  \ \ \dot{y}= y^2(\alpha x + b y) + x^5,
\end{equation}
with $\Lambda = \{ (\alpha, b, c) \in \mathbb{R}^3 : \alpha < 0, c < 0 \}$ belongs to the $\mathcal{G}$-monodromic class and possesses $\eta_1 = 1$. Despite $(1,1) \in W(\mathbf{N}(\mathcal{X}))$, it is shown that in polar coordinates $PV \int_0^{2 \pi} \frac{R_{1}(\varphi)}{G_{0}(\varphi)} d \varphi$ does not exists when $b \neq 0$. Notice that the origin of \eqref{Ej2-Nonlinearity} with parameters in $\{b=0\} \cap \Lambda$ is a reversible center.

Now we try to apply the techniques of this work taking into account that $W(\mathbf{N}(\mathcal{X})) = \{(1,1), (1,2)\}$. First we make the transformation $(x,y) \mapsto (y,-x)$ that reaches the condition $0 \not\in \Omega_{11} = \{\pi/2, 3 \pi/2\}$. Indeed, we obtain $G_0(\varphi) = -c \cos^4\varphi$ and $R_1(\varphi) = b \cos^2\varphi -(c+2\alpha + c \cos(2 \varphi)) \sin(2\varphi)/4$. It is easy to check that there is no $m \geq 2$ such that the function $-F$ of Corollary \ref{teo-eta1-m} has not at most a simple pole at $\pi/2$. Therefore, there is no $m \geq 2$ with $\hat{g}_m$ analytic at $\pi/2$ and so we cannot use Theorem \ref{cor-calculo-eta1} to compute $\eta_1$ using polar coordinates. Indeed, straightforward computations show that $\hat{g}_m$ in $C^\infty$ flat at $\Omega_{11}$.
\newline

Instead, we are going to use $(1,2)$-weighted polar coordinates $(\varphi, \rho)$. In these coordinates we have $G_0(\varphi) = 2 \cos^2\varphi (\cos^4\varphi - \alpha \sin^2\varphi)/(\cos(2\varphi)-3)$ and $R_1(\varphi) = (3 + 8 \alpha + 4 \cos(2\varphi) + \cos(4\varphi)) \sin(2\varphi)/(8 (\cos(2\varphi)-3))$. Now $\Omega_{12} = \{\pi/2, 3 \pi/2\}$ and the function $-F$ of Corollary \ref{teo-eta1-m} has at most a simple pole at $\Omega_{12}$ with ${\rm Res}(-F,\varphi_j^*) = 3-m$. Therefore, for any $m \geq 2$, $\hat{g}_m$ is analytic at $\Omega_{12}$. We continue taking $m=3$ so that $\hat{g}_3(\varphi_j^*) \neq 0$. Indeed we obtain
$$
\hat{g}_3(\varphi) = \frac{(3 - 4 \alpha + 4 (1 + \alpha) \cos(2\varphi) + \cos(4\varphi))^2}{3-\cos(2\varphi)}.
$$

We use now Theorem \ref{teo-Fund-G-jet} in each characteristic direction $\varphi_1^* = \pi/2$ and $\varphi_2^* = 3\pi/2$. We do the blow up $x = \rho u$ with $x=\varphi-\varphi^*_j$ using $\alpha=1$. Then we obtain $R(\varphi_j^* + \rho u, \rho) = \rho^2 Q(u) + o(\rho^{2})$ and $\Theta(\varphi_j^* + \rho u, \rho) = \rho^{2} P(u) + o(\rho^{2})$ for $j=1,2$ obtaining $r=2$ and $Q(u)=(b-u)/2$ and $P(u)=(-2 c + b u - u^2)/2$.

Following Remark \ref{Remar-funcion-a} with $\hat{g}_3(\varphi_j^*) \neq 0$ we take the function $a(\varphi, \rho) = \hat{g}_7(\varphi)$ obtaining $a(\varphi_j^* + \rho u, \rho) = B + O(\rho)$ for $j=1,2$ with $B=16$. Hence $\nu=0$ the second bound of \eqref{eq:positiveDB} holds for any constant $0 < c_B < 16$.

We calculate $D(u) = P(u)- u Q(u) = -c > 0$ on $\Lambda$ and, since it is a constant, we see that the first bound in \eqref{eq:positiveDB} given by $D(u) \ge c_D(1+|u|)^{2}$ is never satisfied by any constant $c_D >0$. In short, we cannot use Theorem \ref{cor-calculo-eta1} to compute $\eta_1$ using $(1,2)$-weighted polar coordinates.

\subsection{An example from \cite{Ga-Ma-Ma} in the class $S_{3\omega}$ revisited}

In the work \cite{Ga-Ma-Ma} it is proved that the origin of family
\begin{equation}\label{ejemplo3-JDE-Armengol}
\dot{x} = b x^2 y + a x y^2 - b y^3 - x^4, \ \ \ \dot{y} = 4 b x y^2 + a y^3 + 2 x^5,
\end{equation}
has monodromic parameter set $\Lambda = \{ (a, b) \in \mathbb{R}^2 : b > 1/4 \}$ and it belongs to the class $S_{3\omega}$. Moreover, in \cite{Ga-Ma-Ma} it is also proved that
\begin{equation}\label{eta1-ARmengol-Skw}
\eta_1 = \exp\left( \frac{2 \pi \sqrt{3}}{3} \frac{a}{b} \right),
\end{equation}
concluding that the origin is an attractor (resp. repeller) focus if $a < 0$ (resp. $a > 0$) while it is a time-reversible center for $a=0$.
\newline

We are going to reproduce \eqref{eta1-ARmengol-Skw} using our method. First we check that $(1,1)$ are weights of the Newton diagram of \eqref{ejemplo3-JDE-Armengol} and that $0 \in \Omega_{11}$, hence changing coordinates $(x,y) \mapsto (y, -x)$ in \eqref{ejemplo3-JDE-Armengol} we get a vector field  $\mathcal{X}$ with $0 \not\in \Omega_{11} = \{\pi/2, 3 \pi/2\}$. Second we see that $\Lambda_{11} = \emptyset$ because the equation $\Theta(\varphi, \rho) = 0$ is polynomial of degree 2 in $\rho$ with discriminant $2 - 64 b + (32 b-2) \cos(2 \varphi) < 0$ for any $\varphi$ in a neighborhood of $\Omega_{11}$ and $b \in \Lambda$.

We get the functions $G_0(\varphi) = b \cos^2\varphi (\cos^2\varphi + 3 \sin^2\varphi)$ and $R_1(\varphi) = a \cos^2\varphi - b \cos\varphi (3 \sin\varphi + \sin(3 \varphi))/2$. Then the function $-F$ of Corollary \ref{teo-eta1-m} with $m \geq 2$ has residue at $\varphi^*=\pi/2$ given by ${\rm Res}(-F,\varphi^*) = (7 - m)/3$ that belongs to $2 \mathbb{N}_0$ if and only if $m = 7$. Up to a nonzero multiplicative constant, the nontrivial solution of \eqref{EDO-GM} with $m = 7$, that is analytic in $[0, 2\pi]$, is
$$
\hat{g}_7(\varphi) = \exp\left( \frac{-6 a}{b} J(\varphi) \right) (\cos^2\varphi + 3 \sin^2\varphi)^5,
$$
where $J(\varphi) = \int_0^{\varphi} (\cos^2\theta + 3 \sin^2\theta)^{-1} \, d \theta$.

\begin{remark}
{\rm We prefer to express $\hat{g}_7$ via the integral $J(\varphi)$ because the primitive $\arctan(\sqrt{3} \tan\varphi)/\sqrt{3}$ of $(\cos^2\theta + 3 \sin^2\theta)^{-1}$ is not continuous in $[0, 2\pi]$ since it has a finite jump discontinuity at $\Omega_{11}$.  We can remove this jump taking the correct branch of the multivalued $\arctan$ function. More precisely we get
$$
\hat{g}_7(\varphi) = (\cos^2\theta + 3 \sin^2\theta)^5 \, \exp\left( \frac{-2 \sqrt{3} \, a}{b} f(\varphi) \right)
$$
with
\begin{equation}\label{ff}
f(\varphi) = \begin{cases}
\arctan(\sqrt{3}\tan\varphi) & \varphi \in [0, \pi/2), \\
\arctan(\sqrt{3}\tan\varphi) + \pi & \varphi \in (\pi/2, 3\pi/2), \\
\arctan(\sqrt{3}\tan\varphi) + 2\pi & \varphi \in (3\pi/2, 2\pi].
\end{cases}
\end{equation}
Notice that we have extended $\hat{g}_m$ at $\Omega_{11}$ by continuity but such extension becomes indeed analytic since the derivative of $\arctan(\sqrt{3}\tan\varphi)$ is the analytic function $\sqrt{3}/(\cos(2 \varphi)-2)$.
}
\end{remark}

We continue applying Theorem \ref{teo-Fund-G-jet} to each $\varphi^* \in \Omega_{11}$. For $\varphi_1^* = \pi/2$ and $\varphi_2^* = 3\pi/2$ we do the blow up $x = \rho u$ with $x=\varphi-\varphi^*_j$, that is, we take $\alpha=1$. Straightforward computations give $R(\varphi_j^* + \rho u, \rho) = \rho^2 Q_j(u) + o(\rho^{2})$ and $\Theta(\varphi_j^* + \rho u, \rho) = \rho^{2} P_j(u) + o(\rho^{2})$ for $j=1,2$ so that $r=2$ where $Q_1(u)=1+b u$, $Q_2(u)=-1+b u$, $P_1(u)=2 - u + 3 b u^2$ and $P_2(u)=2 + u + 3 b u^2$.

According to Remark \ref{Remar-funcion-a}, since $\hat{g}_7(\varphi_j^*) \neq 0$ we take the function $a(x,\rho)$ of Theorem \ref{teo-Fund-G-jet} as $a(x, \rho) = \hat{g}_7(x)$, so that $\nu=0$ and $B_j$ is the constant $B_j(u) = \hat{g}_7(\varphi_j^*) > 0$. In particular, the second bound of \eqref{eq:positiveDB} holds for any constant $0 < c_B < \hat{g}_7(\varphi_j^*)$.

Letting $D_j(u) = P_j(u)- u Q_j(u)$ we will get the first bound in \eqref{eq:positiveDB}, that is, $D_j(u) \ge c_D^{(j)}  \, \beta_j(u)$ with $\beta_j(u) = (1+|u|)^{2}$ for some $c_D^{(j)}>0$. This bound follows using the idea of Remark \ref{Remar-cotas} for the positive functions $D_1(u) = 2 (1 - u + b u^2)$ and $D_2(u) = 2 (1 + u + b u^2)$ by monodromy.

Now the conditions of Theorem \ref{teo-Fund-G-jet} hold and we can use Theorem \ref{cor-calculo-eta1} to compute $\eta_1$ as follows: the functions $\Xi_j(u)$ are $\Xi_1(u)= 6/D_1(u)$ and $\Xi_2(u)= - 6/D_2(u)$, hence $\int_{\mathbb{R}} \Xi_1(u) \, du = - \int_{\mathbb{R}} \Xi_2(u) \, du = 6 \pi/\sqrt{4b-1}$ and $C_7 = 1$ so
$$
\eta_1 = \left( \frac{\hat{g}_m(0)}{\hat{g}_m(2 \pi)} \right)^{\frac{1}{6}} = \exp\left( PV \int_{0}^{2 \pi} \frac{R_1(\varphi)}{G_0(\varphi)} d \varphi   \right) = \exp\left( \frac{2 \pi \sqrt{3}}{3} \frac{a}{b} \right),
$$
yielding \eqref{eta1-ARmengol-Skw}.

\subsection{The Ma\~{n}osa monodromic family I revisited} \label{section-Manosa}

In \cite{Ma} it proved that family
\begin{equation}\label{ejemplo-manyosa}
\dot{x} = x y^2 - y^3 + a x^5, \ \ \dot{y} = 2 x^7 - x^4 y + 4 x y^2 + y^3,
\end{equation}
has a monodromic singularity at the origin for parameters in $\Lambda = \{a \in \mathbb{R} : \Delta(a) := 32 - (1+3 a)^2 >0 \}$. Moreover, there it is shown that
\begin{equation}\label{eta1-manyosa}
\eta_1 = \exp\left( \pi + \frac{4 \pi a}{\sqrt{\Delta(a)}} \right).
\end{equation}
Nevertheless the stability of the origin of \eqref{ejemplo-manyosa} when $a = -31/25$ (corresponding to $\eta_1 = 1$) was left open in \cite{Ma}.
\newline

We will try to apply our method for such family. Since $(1,1)$ are weights of the Newton diagram of \eqref{ejemplo-manyosa} with $0 \in \Omega_{11}$, we first do the change of coordinates $(x,y) \mapsto (y, -x)$ obtaining a vector field  $\mathcal{X}$ with $0 \not\in \Omega_{11} = \{\pi/2, 3 \pi/2\}$. We also claim that $\Lambda_{11} = \emptyset$ because the equation $\Theta(\varphi, \rho) = 0$ is polynomial of degree 4 without odd exponents in $\rho$ with discriminant $-39 + a (2 + a) - (-23 + a (2 + a)) \cos(2 \varphi) < 0$ for any $\varphi$ in a neighborhood of $\Omega_{11}$ and $a \in \Lambda$.

Straightforward computations yield  $G_0(\varphi) = \cos^4\varphi + \sin^2(2 \varphi))$ and $R_1(\varphi) = \cos^2\varphi -3  \cos^3\varphi  \sin\varphi$, hence the function $-F$ of Corollary \ref{teo-eta1-m} has ${\rm Res}(-F,\pi/2) = 2 \in 2 \mathbb{N}_0$. This implies that all the solutions of \eqref{EDO-GM} are analytic at $[0, 2\pi]$ for every $m \in \mathbb{N}$. Actually they are, up to a multiplicative constant,
$$
\hat{g}_m(\varphi) = \exp\left( -2(m-1) J(\varphi) \right) \cos^2\varphi (5- 3 \cos(2 \varphi))^{(m+1)/2},
$$
where $J(\varphi) = \int_0^{\varphi} (5- 3 \cos(2 \theta))^{-1} \, d \theta$.

We do the blow-up of Theorem \ref{teo-Fund-G-jet} to each $\varphi^* \in \Omega_{11}$. For both $\varphi_1^* = \pi/2$ and $\varphi_2^* = 3\pi/2$, the blow up is $x = \rho^2 u$ with $x=\varphi-\varphi^*_j$, hence  $\alpha=2$. We get moreover $R(\varphi_j^* + \rho^2 u, \rho) = \rho^3 Q(u) + o(\rho^{3})$ and $\Theta(\varphi_j^* + \rho^2 u, \rho) = \rho^{4} P(u) + o(\rho^{4})$ for $j=1,2$ so that $r=3$ where $Q(u)=a$, and $P(u) = 2 - (1+a) u + 4 u^2$.

Notice that $\hat{g}_m(\varphi^*) = 0$ and therefore, see  Remark \ref{Remar-funcion-a}, we cannot take the function $a(x,\rho)$ of Theorem \ref{teo-Fund-G-jet} to be simply $\hat{g}_m(x)$. Instead we take $a(x, \rho) = \hat{g}_m(x) + c_* \rho^{\kappa \alpha}$ where $c_* > 0$ is a constant and $\kappa = 2$ is the vanishing multiplicity of $\hat{g}_m$ at $\varphi^* = \pi/2$. In short we take $a(x, \rho) = \hat{g}_m(x) + c_* \rho^{4}$ that is transformed by the blow-up to $a(\varphi_j^* +\rho^2 u, \rho) = \rho^{4} B_j(u) + o(\rho^{4})$ with $\nu=4$ and $B_j(u) = c_* + k_j(m) u^2$ with $k_1(m)= 2^{3 (1 + m)/2} \exp((1 - m) \pi/4)$ and $k_2(m)= 2^{3 (1 + m)/2} \exp(3 (1 - m) \pi/4)$.

The function $D(u) = P(u)- \alpha u Q(u) = 2 - (1+3 a) u + 4 u^2 > 0$ under the monodromic parameters $\Lambda$. We reach the first bound in \eqref{eq:positiveDB} given by $D(u) \ge c_D \, \beta(u)$ with $\beta(u) = (1+|u|^{1/2})^{4}$ for some $c_D > 0$ just using Remark \ref{Remar-cotas} because $D(u)/\beta(u) \to 4$ as $|u| \to \infty$. The second bounds in \eqref{eq:positiveDB} holds and give $B_j(u) \ge c_B \, \beta_j(u)$ with $\beta_j(u) = (1+|u|^{1/2})^4$ for some $c_B > 0$ using again Remark \ref{Remar-cotas} with $B_j(u)/\beta_j(u) \to k_j(m)$ as $|u| \to \infty$.

The hypothesis of Theorem \ref{teo-Fund-G-jet} hold and so we use Theorem \ref{cor-calculo-eta1} to compute $\eta_1$ as follows: the integrals $\int_{\mathbb{R}} \Xi_1(u) \, du = \int_{\mathbb{R}} \Xi_2(u) \, du = (2 a (m-1) \pi)/\sqrt{\Delta(a)}$ and so
$$
\eta_1 = \left( C_m \frac{\hat{g}_m(0)}{\hat{g}_m(2 \pi)} \right)^{\frac{1}{m-1}} = \exp\left( \pi + \frac{4 \pi a}{\sqrt{\Delta(a)}} \right),
$$
yielding \eqref{eta1-manyosa}.

\subsection{A Huelva monodromic family} \label{section-Huelva}

We consider the family of polynomial differential systems
\begin{equation}\label{(a)-Huelva}
\dot{x} = a y^3 + c x y^2 + g x^2 y+ e x^5, \ \ \dot{y} = d y^3 + h x y^2 + f x^4 y+ b x^7,
\end{equation}
corresponding to case (a) of Theorem 13 in \cite{AGR2}, that is with parameter in
\begin{eqnarray}
\Lambda_a &=& \{ a b (g^2+h^2) < 0, (h-g)(h-3 g)>0, (d-c)^2 + 4 a(h-g) < 0, \nonumber \\
 & & (f-3 e)^2 -4 b (h-3 g) < 0\}. \label{La}
\end{eqnarray}
In \cite{AGR2} it is proved that, under these parameter constrains the origin of \eqref{(a)-Huelva} is monodromic, that is $\Lambda_a \subset \Lambda$. We want to mention that the analysis of the monodromy of the origin of \eqref{(a)-Huelva} as well as the characterization of its stability are both quite degenerated problems. The reasons are that \eqref{(a)-Huelva} does not satisfy the sufficient monodromic conditions stated in Theorem 4 of \cite{AGR} neither the non-degeneracy condition of Brunella-Miari type \cite{BM}, see also \cite{GGV}. Additionally it is easy to check that the origin of \eqref{(a)-Huelva} does not belong to the monodromic classes $S_{k\omega}$ of \cite{Ga-Ma-Ma} and class $\mathcal{G}$ of \cite{Ga-Li-Ma-Ma}.

Imposing this monodromic family has no curves of zero angular speed at the origin (parameters in $\Lambda_a \setminus \Lambda_{11}$) and the flow rotates in the positive sense ($h-g > 0$) we determine the stability of the origin at first order in the next result when the resonant condition $h/g \in 2 \mathbb{N} + 1$ holds.

\begin{proposition}
We consider family \eqref{(a)-Huelva} under the parameter restrictions $(\Lambda_a \setminus \Lambda_{11}) \cap \{ h-g > 0 \} \cap \{ 5 \leq h/g \in 2 \mathbb{N} + 1 \}$ where $\Lambda_a \setminus \Lambda_{11} = \Lambda_a \cap \{ S_1 \cup S_2 \}$ where $S_1 = \{ \Delta_1 < 0 \}$ and $S_2 = \{ \Delta_1 = c - d = 0 \}$ being $\Delta_1 = (e - f)^2 + 4 b (g - h)$. Then the Poincar\'e map is $\Pi(\rho) = \eta_1 \rho + o(\rho)$ with
$$
\eta_1 = \exp\left( \left( \frac{2 d g - h (c+d)}{\sqrt{4 a (g - h)-(c - d)^2}} + \frac{f g+e(3 g-2 h)}{\sqrt{4 b (h-3 g)} -(f-3 e)^2} \right)  \, \frac{2 \pi}{g-h} \right).
$$
\end{proposition}

\begin{proof}
This family has Newton diagram with weights $W(\mathbf{N}(\mathcal{X})) = \{ (1,1), (1,3)\}$. In what follows $(\varphi, \rho)$ denotes the polar coordinates. One has $0 \in \Omega_{11}$, so we do the change of coordinates $(x,y) \mapsto (y, -x)$ and the resulting vector field  $\mathcal{X}$ has $0 \not\in \Omega_{11}$. The polar vector field $\mathcal{Z}$ of the transformed field is $\mathcal{Z} = (G_0(\varphi) + O(\rho)) \partial_\varphi + (R_1(\varphi) \rho + O(\rho^2)) \partial_\rho$ with $G_0(\varphi) = \cos^2\varphi \, \hat{G}_0(\varphi)$ with $\hat{G}_0(\varphi) = -a \cos^2\varphi + (c-d) \cos\varphi \sin\varphi + (h-g) \sin^2\varphi$ and $R_1(\varphi) = \cos\varphi (d \cos^3\varphi  - (a+h) \cos^2\varphi \sin\varphi + c \cos\varphi \sin^2\varphi - g \sin^3\varphi)$. It is easy to check that $\hat{G}_0$ has no real roots when restricted to $\Lambda_s$, hence $\Omega_{11} = \{\pi/2, 3 \pi/2\}$. Actually it follows that $G_0 \geq 0$ if and only if $h-g > 0$ and we restrict the following analysis to this choice of the flow sense rotation without loss of generality.

The function $-F$ of Corollary \ref{teo-eta1-m} has residue ${\rm Res}(-F,\pi/2) = (g - 2 h + g m)/(g - h)$. For simplicity we take the residue vanishes, hence $m =-1 + 2h/g$ and we impose the resonance condition $h/g \geq 4$ be an integer compatible with $\Lambda_a$. With this choice the non-trivial solutions $\hat{g}_m$ of \eqref{EDO-GM} are analytic on $[0, 2\pi]$ and $\hat{g}_m|_{\Omega_{11}} \neq 0$ so $\kappa = 0$. Actually they are, up to a multiplicative constant, $\hat{g}_m(\varphi) = \exp(J(\varphi)) K^{1+h/g}(\varphi)$, where
\begin{equation}\label{J-Huelva}
J(\varphi) = \frac{2 (-2 d g+ (c + d) h)}{g}  \int_{0}^{\varphi} K^{-1}(\theta) d \theta,
\end{equation}
and $K(\varphi) = a + g-h + (a -g+h) \cos(2 \varphi) + (d-c) \sin(2 \varphi) < 0$ since $\hat{G}_0(\varphi)/K(\varphi) = -1/2$.

We are going to prove that the parameter set $\Lambda_a\setminus \Lambda_{11}$ where there are no curves of zero angular speed for the origin of \eqref{(a)-Huelva} is characterized by $\Lambda_a \setminus \Lambda_{11} = \Lambda_a \cap \{ S_1 \cup S_2 \}$ where $S_1 = \{ \Delta_1 < 0 \}$ and $S_2 = \{ \Delta_1 = c - d = 0 \}$ being $\Delta_1 = (e - f)^2 + 4 b (g - h)$. The proof relies in the fact that the equation $\Theta(\varphi, \rho) = 0$ is a polynomial in $\rho$ of degree 4 with even exponents and has a discriminant $\delta(\varphi) = \cos^2\varphi \sin^8\varphi (4 a b \cos^2\varphi + 4 b (d -c) \sin\varphi \cos\varphi + \Delta_1 \sin^2\varphi)$. From the Taylor expansion at $\pi/2$ of $\delta$ we have $\delta(\theta) = \Delta_1  (\varphi-\pi/2)^2 + \cdots$ and therefore, in a neighborhood of $\varphi = \pi/2$, we have $\delta(\varphi) < 0$ when  $\Delta_1 < 0$. The degenerate case $\Delta_1 = 0$ with $c=d$ has not associated curves of zero angular speed since in this case $\delta(\theta) = 4 a b (\varphi-\pi/2)^4 + \cdots$ and $a b < 0$ on $\Lambda_a$ and so $\delta(\theta) < 0$ too. In what follows we restrict the parameters to
\begin{equation}\label{mon-rest.ej}
(\Lambda_a \setminus \Lambda_{11}) \cap \{ h-g > 0 \} \cap \{ 4 \leq h/g \in \mathbb{N} \}.
\end{equation}

The blow-up of Theorem \ref{teo-Fund-G-jet} to each characteristic direction $\varphi_1^* = \pi/2$ and $\varphi_2^* = 3\pi/2$ is $x = \rho^2 u$ with $x=\varphi-\varphi^*_j$, hence  $\alpha=2$. The outcome is $R(\varphi_j^* + \rho u, \rho) = \rho^3 Q(u) + o(\rho^{3})$ and $\Theta(\varphi_j^* + \rho u, \rho) = \rho^{4} P(u) + o(\rho^{4})$ for $j=1,2$ so that $r=3$ where $Q(u)=e + g u$, and $P(u) =b + (f-e) u+ (h-g) u^2$.

The function $D(u) = P(u)- \alpha u Q(u) = b + (f-3e) u+ (h-3g) u^2 > 0$ when the parameters lie in \eqref{mon-rest.ej}. The first bound in \eqref{eq:positiveDB}, that is $D(u) \ge c_D \, \beta(u)$ with $\beta(u) = (1+\sqrt{|u|})^{4}$ for some $c_D > 0$, follows by Remark \ref{Remar-cotas} because $D(u)/\beta(u) \to h-3 g > 0$ as $|u| \to \infty$.

Since $\hat{g}_m|_{\Omega_{11}} \neq 0$, according to Remark \ref{Remar-funcion-a}, we take $a(\varphi) = \hat{g}_m(\varphi)$. Then, following Theorem \ref{teo-Fund-G-jet}, we have $a(\varphi_j^* +\rho u, \rho) = B_j + O(\rho)$ with $B_j=\hat{g}_m(\varphi_j^*)$ constant and $\nu=0$. In particular, the second bound of \eqref{eq:positiveDB} holds for any constant $0 < c_B < \hat{g}_m(\varphi_j^*)$. More specifically one has $B_j= \exp(J(\varphi_j^*)) (2 (g-h))^{1+h/g}$. In order that $B_j >0$, and recalling that $h-g > 0$ and $h/g \geq 4$ is in $\mathbb{Z}$, we need to impose that the exponent $1+h/g$ be in $2 \mathbb{N}$. In summary we need $h/g \geq 5$ and odd and then we replace the set $\{ 4 \leq h/g \in \mathbb{N} \}$ in \eqref{mon-rest.ej} by $\{ 5 \leq h/g \in 2 \mathbb{N} + 1 \}$

The concrete values of $B_j$ can be explicitly computed via the involved quadratures \eqref{J-Huelva} after the change of variables $t= \tan\theta$ and noticing that the discriminant $\Delta_2$ of the polynomial $a +(d-c) t +(g-h) t^2$ in the next denominator is  $\Delta_2 = (d-c)^2 - 4 a(g-h) < 0$ on $\Lambda_a$. The outcome is
\begin{eqnarray*}
J(\pi/2) &=& \frac{-2 d g+ (c + d) h}{g} \int_{0}^{\infty} \frac{d t}{a +(d-c) t +(g-h) t^2} \\
 &=& - 2 \frac{-2 d g+ (c + d) h}{g \sqrt{-\Delta_2}} \left( \frac{\pi}{2} + \arctan\left( \frac{d-c}{\sqrt{-\Delta_2}} \right) \right), \\
J(3\pi/2) &=& - 2 \frac{-2 d g+ (c + d) h}{g \sqrt{-\Delta_2}} \left( \frac{3 \pi}{2} + \arctan\left( \frac{d-c}{\sqrt{-\Delta_2}} \right) \right).
\end{eqnarray*}

Now we are in position to apply Theorem \ref{teo-Fund-G-jet} and also Theorem \ref{cor-calculo-eta1} to compute $\eta_1$. We compute
$$
\Xi_1(u) = \Xi_2(u) = - \frac{3 e g + f g - 2 e h}{g} \, \frac{1}{b + (f- 3 e) u + (h-3 g) u^2}
$$
that is a bounded function in $\mathbb{R}$. Therefore,
$$
\int_{\mathbb{R}} \Xi_1(u) du = - \frac{3 e g + f g - 2 e h}{g} \frac{2 \pi}{\sqrt{\Delta_3}}
$$
where $\Delta_3 = 4 b (h-3 g) - (f-3 e)^2 > 0$ in $\Lambda_a$. Then $C_m = \exp\left(\int_{\mathbb{R}} \Xi_1(u) du \right)^2$ and $\hat{g}_m(0)/\hat{g}_m(2 \pi) = \exp(-J(2 \pi))$ with
$$
J(2\pi) = - 4 \pi \frac{-2 d g+ (c + d) h}{g \sqrt{-\Delta_2}}.
$$
In summary we get from
$$
\eta_1 = \left( C_m \frac{\hat{g}_m(0)}{\hat{g}_m(2 \pi)} \right)^{\frac{1}{m-1}}
$$
the value of $\eta_1$ given in the proposition. This finishes the proof.
\end{proof}

\section{Proofs}

In this section we will prove the main results of this work.

\subsection{Proof of Theorem \ref{teo-Fund-GPi2}}

First we need a technical lemma.

\begin{lemma}\label{structure-G=0}
In the restricted monodromic parameter subset $\Lambda \setminus \Lambda_{pq}$, any solution $\mathcal{G}^*(\varphi, \rho)$ of \eqref{edp-G*} analytic in $\mathcal{R} \setminus \Theta^{-1}(0)$ satisfies $(\mathcal{G}^*)^{-1}(0) \cap \mathcal{R} = \emptyset$ with $I$ sufficiently small.
\end{lemma}
\begin{proof}
We see from \eqref{edp-G*} that the set $(\mathcal{G}^*)^{-1}(0)$ must be invariant for the flow of $\mathcal{Z}^*$. Therefore, either $(\mathcal{G}^*)^{-1}(0)$ contains spirals tending to $\rho = 0$ as $\varphi$ goes to $\infty$ or $-\infty$ (in the focus case) or it contains a sequence of periodic orbits of $\mathcal{Z}^*$ accumulating at $\rho=0$ (in the center case). In the first case the lemma follows as in Remark \ref{Rem-G-cero-linear} taking into account that $0 \leq \varphi \leq 2 \pi$ in $\mathcal{R}$. The second situation is not possible since it would imply that $(\mathcal{G}^*)^{-1}(0) \cap \mathcal{R}$ has a sequence of components accumulating at $\rho=0$ and, taking an fixed angle $\hat\varphi \in [0, 2 \pi]$, the function $\mathcal{G}^*(\hat\varphi, \rho)$ will be real analytic at $\rho=0$ and would have a sequence of zeros accumulating at $\rho=0$, a contradiction.
\end{proof}

\begin{proof}[Proof of Theorem \ref{teo-Fund-GPi2}]
The functions $\mathcal{G}$ and $\mathcal{G}^*$ are related by $\mathcal{G}^* = \mathcal{G}/\Theta$.

Take $0 \not\in \Omega_{pq}$ without loss of generality and assume an initial condition $\mathcal{G}^*(0, \rho)$ analytic in $I$. We notice that the vector field $\mathcal{Z}^*$ is analytic in $\mathcal{R} \setminus \Theta^{-1}(0)$ hence the linear first order partial differential equation \eqref{edp-G*} is also analytic in $\mathcal{R} \setminus \Theta^{-1}(0)$. Since $\{\varphi = 0 \} \subset \mathcal{R}$ is a transversal to the flow of $\mathcal{Z}^*$, we use the characteristics method to construct the unique analytic solution $\mathcal{G}^*(\varphi, \rho)$ of \eqref{edp-G*} with given initial condition $\mathcal{G}^*(0, \rho) = f(\rho)$ analytic in $I$. Clearly, by monodromy, this $\mathcal{G}^*$ is defined in the hole $\mathcal{R} \setminus \{ \rho=0 \}$.

We claim that $\mathcal{G}^*$ can be uniquely analytically extended to $\rho=0$ except at the singularities $(\varphi, \rho) = (\varphi^*, 0)$ for all $\varphi^* \in \Omega_{pq}$. More precisely, since $0 \not\in \Omega_{pq}$, by Cauchy-Kovalevskaya Theorem, there is a unique local analytic solution $\hat{\mathcal{G}}(\varphi, \rho)$, defined in a neighborhood $U$ of the non-singular point $(\varphi, \rho) = (0, 0)$, of the Cauchy problem formed with  \eqref{edp-G*} and any initial condition $\hat{\mathcal{G}}(0, \rho)$ analytic at $\rho=0$. If we now take $\hat{\mathcal{G}}(0, \rho) = f(\rho)$, by the unicity of the solutions $\mathcal{G}^*$ and $\hat{\mathcal{G}}$ of \eqref{edp-G*} it follows that the restriction $\mathcal{G}^*|_{U \cap \mathcal{R}} = \hat{\mathcal{G}}|_{U \cap \mathcal{R}}$. Indeed, $\hat{\mathcal{G}}$ can be analytically extended along the segment $L = \{ (\varphi, \rho) = (\varphi, 0) : 0 \leq \varphi < \varphi_1^* \}$ being $\varphi_1^* = \min\{ \Omega_{pq} \}$. In this way we get that $\hat{\mathcal{G}}$ is the analytic extension of $\mathcal{G}^*$ to $L \subset \{\rho=0\}$. It may happen that the radius of convergence of $\hat{\mathcal{G}}(\varphi, \rho)$ goes to $0$ as $\varphi \to \varphi_1^*$ from the left. Repeating this analytic extension argument to the connected set $\mathcal{R} \setminus \Theta^{-1}(0)$ we prove that the real analytic function $\mathcal{G}^*(\varphi, \rho)$ is defined on $\mathcal{R} \setminus \Theta^{-1}(0)$ and therefore the same holds for the function $\mathcal{G}(\varphi, \rho)$.
\newline

Now we are going to prove the second part of the theorem. First we focus on the zero set $(\mathcal{G}^*)^{-1}(0)$. The differential equation for $g_1^*$ implies that it is a constant $g_1^*(\varphi) = C_1$ that we chose, without loss of generality, $C_1 \neq 0$. Therefore, for any $\varphi \not\in \Omega_{pq}$, one has $\mathcal{G}^*(\varphi, \rho) = \rho \, U^*(\varphi, \rho)$ where $U^*(\varphi, 0) \neq 0$. But this argument does not prevent the existence of a component of $(\mathcal{G}^*)^{-1}(0)$ emerging from the point $(\varphi^*, 0)$. To remove this possibility we apply Lemma \ref{structure-G=0}.

Let $\rho_0 > 0$ and $\Delta \rho_0 > 0$ be sufficiently small, and let $\Psi(\varphi; \rho_0)$ be the solution of equation
\begin{equation}
\frac{d \rho}{d\varphi} \, = \, \mathcal{F}(\varphi, \rho), \label{eq3}
\end{equation}
with initial condition $\Psi(0; \rho_0) = \rho_0$. We consider a simply connected domain $D \subset \mathcal{R}$ whose boundary $\partial D$ is the union of the curves $\Gamma_{\rho_0} = \{ (\varphi, \Psi(\varphi; \rho_0)) : 0 \leq \varphi \leq 2 \pi \}$,
$\Gamma_{\rho_0 + \Delta \rho_0} = \{ (\varphi, \Psi(\varphi; \rho_0 + \Delta \rho_0)) : 0 \leq \varphi \leq 2 \pi \}$,
and the vertical segments $\ell_0 = \{ (0, \rho) : \rho_0 \leq r \leq \rho_0 + \Delta \rho_0 \}$ and
$\ell_{2 \pi} = \{ (2 \pi, \rho) : \Pi(\rho_0) \leq \rho \leq \Pi(\rho_0 + \Delta \rho_0) \}$.

The condition ${\rm div}\left( \mathcal{Z}^*/ \mathcal{G}^* \right) \equiv 0$ holds in $D$ because $\Theta^{-1}(0) \cup (\mathcal{G}^*)^{-1}(0)$ is empty in $D$. Taking the positive counterclockwise orientation on $\partial D$, we apply Green's theorem in the domain $D$ to the rescaled vector field $\mathcal{Z}^*/ \mathcal{G}^*$ whose associated differential 1--form is $\Omega = (\mathcal{F} d \varphi - d \rho)/ \mathcal{G}^*$:
\begin{equation}\label{cent-5}
\int\int_D {\rm div}\left( \frac{\mathcal{Z}^*}{\mathcal{G}^*} \right) \, d \rho d \varphi = \int_{\Gamma_{\rho_0}}  \Omega + \int_{\ell_{2 \pi}} \Omega + \int_{\Gamma_{\rho_0+ \Delta \rho_0}}  \Omega + \int_{\ell_0} \Omega.
\end{equation}
Clearly, $\int_{\Gamma_{\rho_0}}  \Omega = \int_{\Gamma_{\rho_0+ \Delta \rho_0}}  \Omega = 0$ since $\Omega|_{\Gamma_{\rho_0}} = \Omega|_{\Gamma_{\rho_0+ \Delta \rho_0}} \equiv 0$. Moreover,
$$
\int_{\ell_0} \Omega = - \int_{\rho_0 + \Delta \rho_0}^{\rho_0} \frac{d \rho}{\mathcal{G}^*(0, \rho)},
\qquad
\int_{\ell_{2 \pi}} \Omega = - \int_{\Pi(\rho_0)}^{\Pi(\rho_0+ \Delta \rho_0)} \frac{d \rho}{\mathcal{G}^*(2 \pi, \rho)}.
$$
Equation \eqref{cent-5} becomes
\begin{eqnarray}\label{cent-5*}
\int_0^{2 \pi} d \varphi \int_{\Psi(\varphi; \rho_0)}^{\Psi(\varphi; \rho_0 + \Delta \rho_0)} {\rm div}\left( \frac{\mathcal{Z}^*}{\mathcal{G}^*} \right) \, d \rho  &=& \mathcal{P}_0(\rho_0+ \Delta \rho_0) - \mathcal{P}_0(\rho_0) + \\
 & & \mathcal{P}_{2 \pi}(\Pi(\rho_0)) - \mathcal{P}_{2 \pi}(\Pi(\rho_0+ \Delta \rho_0)), \nonumber
\end{eqnarray}
where $\mathcal{P}_\theta(\rho)$ is a primitive of $1/\mathcal{G}^*(\theta, \rho)$ for a fixed angle $\theta$.

We divide \eqref{cent-5*} by $\Delta \rho_0$ and consider the limit $\Delta \rho_0 \to 0$. Regarding the right-hand side of \eqref{cent-5*}, we observe that
\begin{eqnarray}
\lim_{\Delta \rho_0 \to 0} \frac{\mathcal{P}_0(\rho_0+ \Delta \rho_0)-\mathcal{P}_0(\rho_0)}{\Delta \rho_0} &=& \mathcal{P}_0'(\rho_0) = \frac{1}{\mathcal{G}^*(0, \rho_0)}, \label{limitdr0} \\
\lim_{\Delta \rho_0 \to 0} \frac{\mathcal{P}_{2 \pi}(\Pi(\rho_0)) - \mathcal{P}_{2 \pi}(\Pi(\rho_0+ \Delta \rho_0))}{\Delta \rho_0} &=& - \frac{d}{d \rho} \left(\mathcal{P}_{2 \pi} \circ \Pi \right)(\rho_0) \nonumber \\
  &=& - \mathcal{P}_{2 \pi}'(\Pi(\rho_0)) \Pi'(\rho_0) = -\frac{\Pi'(\rho_0)}{\mathcal{G}^*(2 \pi, \Pi(\rho_0))}. \nonumber
\end{eqnarray}
Therefore we obtain
$$
\lim_{\Delta \rho_0 \to 0} \frac{1}{\Delta \rho_0} \int_0^{2 \pi} d \varphi \int_{\Psi(\varphi; \rho_0)}^{\Psi(\varphi; \rho_0 + \Delta \rho_0)} {\rm div}\left( \frac{\mathcal{Z}^*}{\mathcal{G}^*} \right) \, d \rho = \frac{1}{\mathcal{G}^*(0, \rho_0)} -\frac{\Pi'(\rho_0)}{\mathcal{G}^*(2 \pi, \Pi(\rho_0))}.
$$
The left-hand side of the former equation is $0$ because ${\rm div}\left( \mathcal{Z}^*/\mathcal{G}^*\right) \equiv 0$ in $D$ proving equation
\begin{equation}\label{G-fundamental}
\mathcal{G}^*(2 \pi, \Pi(\rho_0)) = \mathcal{G}^*(0, \rho_0) \, \Pi'(\rho_0).
\end{equation}
Since $\mathcal{G}^* = \mathcal{G}/\Theta$ we deduce \eqref{G-fundamental2}.
\end{proof}

\subsection{Proof of Theorem \ref{teo-Fund-G-jet}}

First we need a technical lemma that allows us to obtain lower bounds of real analytic functions in two variables from the bounds of its leading quasi-homogenous terms.

\begin{lemma} \label{lem:cotas-principales}
Let $f(x,\rho)$ be real analytic near $(0,0)$, and let $f(x,\rho) = f_d(x,\rho) + \cdots$ its $(\alpha,1)$-quasihomogeneous expansion, that is $f_d(x,\rho)=\sum_{\alpha i+j=d} f_{ij}x^i\rho^j$ and $d\ge0$. We define the compact set
\[
 \mathcal{K}=\{(X,R)\in\mathbb R\times[0,\infty): |X|^{1/\alpha}+R=1\}.
\]
Suppose there is a constant $c_0>0$ such that $f_d(X,R)\ge c_0$ for every $(X,R)\in \mathcal{K}$. Then there exists $\varepsilon>0$ such that
\[
 f(x,\rho)\ge \frac{c_0}{2}
 \bigl(|x|^{1/\alpha}+\rho\bigr)^d
\]
whenever $\rho\ge0$ and $0<|x|^{1/\alpha}+\rho<\varepsilon$.
\end{lemma}

\begin{proof}
For $(x,\rho)\ne(0,0)$ with $\rho\ge0$, we define $w := |x|^{1/\alpha}+\rho$ and we take the normalization
\[
 X=\frac{x}{w^\alpha},\qquad R=\frac{\rho}{w}.
\]
This normalization satisfies $(X,R)\in \mathcal{K}$ and $(x,\rho)=(w^\alpha X,wR)$. Writing the convergent series for $f(x, \rho) = \sum_{i,j}f_{ij} x^i \rho^j$, we have
\[
 w^{-d}f(w^\alpha X,wR)
 =\sum_{i,j}f_{ij}w^{\alpha i+j-d}X^iR^j.
\]
No negative powers of $w$ occur. The terms with $\alpha i+j=d$ sum to $f_d(X,R)$; all remaining terms tend uniformly to zero on $\mathcal{K}$ as $w$ tends to zero, that is,
\begin{equation}\label{eq:uniform-normalized-remainder}
 \lim_{t \to 0^+} \sup_{(X,R)\in \mathcal{K}}  \left|t^{-d}f(t^\alpha X,tR)-f_d(X,R)\right|=0.
\end{equation}
Therefore, for sufficiently small $w$, it follows that $\left|w^{-d}f(w^\alpha X, w R)- f_d(X,R)\right|$ is small enough, hence
\[
 w^{-d}f(w^\alpha X,wR)\ge f_d(X,R)-\frac{c_0}{2} \ge\frac{c_0}{2},
\]
where in the last step we have used the hypothesis $f_d(X,R) \ge c_0$. The proof finishes multiplying by $w^d$.
\end{proof}

\begin{proof}[Proof of Theorem \ref{teo-Fund-G-jet}]
That the restriction of $\mathcal{G}$ to $I_{pq} = [0, 2 \pi] \setminus \Omega_{pq}$ admits the convergent Taylor expansion \eqref{Taylor-G} is just as a consequence of the analyticity of $\mathcal{G}$ on $\mathcal{R} \setminus \Theta^{-1}(0)$ already proved in Theorem \ref{teo-Fund-GPi2}. On the other hand, inserting \eqref{Taylor-G} into the partial differential equation \eqref{edp-G} shows that both $g_m$ and $\hat{g}_m$ are solutions of the same linear homogeneous singular differential equation \eqref{EDO-GM}, hence $g_m(\varphi)/\hat{g}_m(\varphi)$ is a non-zero constant in each connected component of $I_{pq}$. The former makes sense to the definitions of the constants $C_{\pm}$.
\newline

Notice that, since $\alpha > 0$ there is always a geometric overlapping region where $\rho \to 0$,  $|u| \to \infty$ and $\rho^\alpha |u| \to 0$.
\newline

Let now $Z =\hat{\Theta}(x, \rho) \partial_x+ \hat{R}(x, \rho) \partial_\rho$ the vector field analytic near $(x,\rho)=(0,0)$ obtained by the translation $x = \varphi-\varphi^*$ of $\mathcal{Z}$. Let $a(x,\rho)$ be analytic in a neighborhood of the origin and define $A(x,\rho)=\rho^m a(x,\rho)$ and $E = (Z(A)-A \, {\rm div} Z) / A$ in the region where $A\ne0$. We define $V$ by $\hat{R} = \rho V$, hence ${\rm div} Z = \hat\Theta_x+V+\rho V_\rho$. A direct calculation gives
\begin{equation}\label{eq:residue}
 E=\frac{N}{a},\qquad
 N=\hat\Theta a_x+\rho V a_\rho+(m-1)V a - a \hat\Theta_x-a\rho V_\rho.
\end{equation}
We take $a(x, \rho)$ such that $a(x,0) = \hat{g}_m(\varphi^*+x)$. Evaluating $N$ at $\rho=0$ gives $G_0(\varphi^*+x) \hat{g}_m'(\varphi^*+x) + ((m-1) R_1(\varphi^*+x) - G_0'(\varphi^*+x)) \hat{g}_m(\varphi^*+x)$. This expression identically vanishes since $\hat{g}_m$ satisfies \eqref{EDO-GM} so
\begin{equation}\label{eq:cancel}
 N(x,0)=0.
\end{equation}

Where $A \neq 0$, from the definition of $E$ and using the equation $Z(G) = G \, {\rm div}(Z)$ that satisfies $G(x,\rho) = \mathcal{G}(\varphi^*+x, \rho)$ we get the identity
\begin{equation}\label{eq:logratio}
 Z\!\left(\log\left|\frac{G}{A}\right|\right)=-E.
\end{equation}
We use the logarithm of the absolute value so as not to require $G$ to be positive since the sign of a nonzero solution of the homogeneous linear partial differential equation is preserved along each trajectory of $Z$.

Let $\Psi(t; \varepsilon, \rho_0) = (x(t; \varepsilon, \rho_0), \rho(t; \varepsilon, \rho_0))$ be the flow of $Z$ with initial condition $\Psi(0; \varepsilon, \rho_0) = (- \varepsilon, \rho_0)$.
Let $\gamma_{\rho_0}$ be an arc of orbit of $Z$ from $\Sigma_{-}^\varepsilon$ to $\Sigma_{+}^\varepsilon$ connecting the transversal sections $\Sigma_{\pm}^\varepsilon = \{ x = \pm \varepsilon \}$. More specifically, $\gamma_{\rho_0} = \{ \Psi(t; \varepsilon, \rho_0) : 0 \leq t \leq \tau \}$ with $\Psi(0; \varepsilon, \rho_0) \in \Sigma_{-}$ and $\Psi(\tau; \varepsilon, \rho_0) \in \Sigma_{+}$, with $\varepsilon > 0$ and $\rho_0 > 0$ small enough and $\tau = \tau(\varepsilon, \rho_0)$. Of course, by monodromy, $\lim_{\rho_0 \to 0^+} \max_{0 \leq t \leq \tau} \{ \rho(t; \varepsilon, \rho_0) \} = 0$.

On the other hand, we will be interested in the integration \eqref{eq:logratio} along $\gamma_{\rho_0}$ parameterized by $t$. Of course $G/A$ is well defined on $\gamma_{\rho_0}$ when $A \neq 0$. Performing this integration yields
\begin{equation}\label{eq:connectionfinite}
 \log\frac{|G/A|^{+}(\varepsilon, \rho_0)}{|G/A|^{-}(\varepsilon, \rho_0)} = - \int_{\gamma_{\rho_0}} E\,  dt := - \int_{0}^{\tau(\varepsilon, \rho_0)} E \circ \Psi(t; \varepsilon, \rho_0) \, dt.
\end{equation}
where we have defined $|G/A|^{-}(\varepsilon, \rho_0) = |G/A| \circ \Psi(0; \varepsilon, \rho_0)$ and $|G/A|^{+}(\varepsilon, \rho_0) = |G/A| \circ \Psi(\tau; \varepsilon, \rho_0)$. Since
\begin{eqnarray*}
\lim_{\rho_0 \to 0^+} \frac{G}{A}(- \varepsilon, \rho_0) &=& \lim_{\rho_0 \to 0^+} \frac{g_m(\varphi^* - \varepsilon) \rho_0^m + O(\rho_0^{m+1})}{a(\varphi - \varepsilon, \rho_0) \rho_0^m} = \frac{g_m(\varphi^* - \varepsilon)}{\hat{g}_m(\varphi^* - \varepsilon)} = C_{-}, \\
\lim_{\rho_0 \to 0^+} \frac{G}{A}(\varepsilon, \rho(\tau; \varepsilon, \rho_0)) &=& \lim_{\rho_0 \to 0^+} \frac{g_m(\varphi^* + \varepsilon) (\rho(\tau; \varepsilon, \rho_0))^m + O((\rho(\tau; \varepsilon, \rho_0))^{m+1})}{a(\varphi + \varepsilon, \rho_0) (\rho(\tau; \varepsilon, \rho_0))^m} \\
 &=& \frac{g_m(\varphi^* + \varepsilon)}{\hat{g}_m(\varphi^* + \varepsilon)} = C_{+},
\end{eqnarray*}
where $C_{+} C_{-} > 0$ since the function $H(t) = G/A \circ \Psi(t; \varepsilon, \rho_0)$ satisfies $H(\tau) = H(0) \exp\left(\int_0^\tau - E dt\right)$ and so $H(\tau)$ has the same sign than $H(0)$. It follows that the limit when $\rho_0 \to 0^+$ in \eqref{eq:connectionfinite} exists and gives
\begin{equation}\label{ec-LogH-C}
\log \frac{C_-}{C_{+}}  = \lim_{\rho_0 \to 0^+} \int_{\gamma_{\rho_0}} E \, dt.
\end{equation}
The issue now is to compute the former limit. To this end we will use dominate convergence: we will perform the blow-up $x = \rho^\alpha u$ that allows to parameterize by $u$ the orbits $\gamma_{\rho_0}$, after we will prove the pointwise convergence of the integrand and finally we will find an integrable majorant independent of $\rho_0$.

By hypothesis, $\hat{R}(\rho^\alpha u, \rho) = Q(u) \rho^r  + O(\rho^{r+1})$ and $\hat{\Theta}(\rho^\alpha u, \rho) = P(u) \rho^{r+\alpha-1}  + O(\rho^{r+\alpha}))$ and $a(\rho^\alpha u, \rho) = B(u) \rho^\nu  + O(\rho^{\nu+1})$ with $\nu \geq 0$. On the other hand, derivating respect to the time $t$ the relation $x=\rho^\alpha u$ and using the vector field $\dot{x} = \hat{\Theta}(x, \rho)$, $\dot{\rho} = \hat{R}(x, \rho)$ we deduce that $\dot{u} = \hat{\Theta}/\rho^\alpha - \alpha \hat{R} u/\rho = T(x, \rho)/\rho^\alpha$ where $T:= \hat{\Theta}-\alpha \rho^{\alpha-1} \hat{R} u = \hat{\Theta}-\alpha x V$. If $T \neq 0$ along $\gamma_{\rho_0}$ then $dt = \rho^\alpha/T(x,\rho) du$. It follows that
\begin{eqnarray}
\int_{\gamma_{\rho_0}} E\,  dt &=& \int_{J_{\rho_0}} \frac{\rho^\alpha E}{T} \,  du  \nonumber \\
   &=& \int_{J_{\rho_0}} \frac{\rho^\alpha(u; \rho_0) N(\rho^\alpha(u; \rho_0) u, \rho(u; \rho_0))}{a(\rho^\alpha(u; \rho_0) u, \rho(u; \rho_0)) \, T(\rho^\alpha(u; \rho_0) u, \rho(u; \rho_0))} \,  du \label{int-gamma}
\end{eqnarray}
where $J_{\rho_0}$ is the interval of $u$-values that $\gamma_{\rho_0}$ crosses and we have parameterized the $\rho$ on $\gamma_{\rho_0}$ by $u$ as $\rho(u; \rho_0)$. Notice that $J_{\rho_0} \to \mathbb{R}$ when $\rho_0 \to 0^+$ because the geometry of the blow-up implies that when $\rho \to 0$ and $\rho^\alpha u$ becomes finite then $|u| \to \infty$.

We are going to apply Lemma \ref{lem:cotas-principales} to the analytic functions $T$ and $a$. We set $w=|x|^{1/\alpha}+\rho=\rho(1+|u|^{1/\alpha})$. From the bounds \eqref{eq:positiveDB} of the leading terms, after shrinking the neighborhood if necessary, we get
\begin{equation}\label{eq:lower}
T(x,\rho) \ge c_T w^h, \qquad a(x,\rho)\ge c_a w^\nu
\end{equation}
for positive constants $c_T$ and $c_a$. In particular, $u$ is strictly increasing since $\dot{u} = T/\rho^\alpha > 0$ on $\gamma_{\rho_0}$ and so $J_{\rho_0} \to (-\infty, \infty)$ when $\rho_0 \to 0$. Moreover, $A=\rho^m a$ does not vanish during the transition on $\gamma_{\rho_0}$.

By \eqref{eq:cancel}, each monomial $n_{ij}x^i\rho^j$ of $N$ satisfies $j\ge1$. Moreover, using the $(\alpha, 1)$-quasihomogeneous degrees of the leading terms in $\hat{R}$, $\hat{\Theta}$ and $a$ that are $r$, $r+\alpha-1$ and $\nu$, respectively, we can check in the definition of $N$ given in \eqref{eq:residue} that all the terms in $N$ have $(\alpha, 1)$-quasihomogeneous degree at least $K=r+\nu-1$, that is $\alpha i+j\ge K$. Regarding the integrand  of \eqref{int-gamma}, for a monomial of $N$, using \eqref{eq:lower}, $|x|\le w^\alpha$, and $\rho/w=v^{-1}$ with $v:= 1+|u|^{1/\alpha}$, we obtain
\[
 \frac{\rho^\alpha |x|^i\rho^j}{aT} \leq \frac{w^{\alpha+j} v^{-(\alpha+j)} w^{\alpha i}}{c_a w^\nu c_T w^{r+\alpha-1}} = K_*^{-1} w^{\alpha i+j-K} v^{-(\alpha+j)} \leq K_*^{-1} w^{\alpha i+j-K} v^{-(\alpha+1)},
\]
with $K_* := c_a c_T$, where in the last step we have used that $j\ge1$ and $v \geq 1$. Since $\alpha i+j-K \ge 0$ and $w$ is uniformly bounded on $\gamma_{\rho_0}$, hence we only consider the region $0 \leq w \leq \delta$, we have $w^{\alpha i+j-K} \leq \delta^{\alpha i+j-K}$. Therefore
\begin{equation}\label{eq:majorant}
\left|\frac{\rho^\alpha N}{aT}\right| \leq  K_*^{-1} v^{-(\alpha+1)} \sum_{i,j} |n_{ij}| \delta^{\alpha i+j-K} = C v^{-(\alpha+1)} = C (1+|u|^{1/\alpha})^{-\alpha-1},
\end{equation}
where we have defined $C = K_*^{-1} \delta^{-K} \sum_{i,j} |n_{ij}| \delta^{\alpha i+j}$. We remark that the series that defines $C$ converges because the analytic series for $N$ is absolutely convergent in a neighborhood of the point $(x, \rho) = (\delta^\alpha, \delta)$.

We emphasize that the majorant \eqref{eq:majorant} is an integrable function on $\mathbb{R}$ because has the asymptotic $1/|u|^{1+1/\alpha}$ when $|u| \to \infty$ with exponent greater that 1.

On the other hand, the function $N(x, \rho)$ is a sum of functions, see \eqref{eq:residue}, all of them having leading $(\alpha,1)$-quasihomogeneous terms of degree $r +\nu -1$. Straight calculations with the terms in $N$ give: $a_x(\rho^\alpha u, \rho) = B'(u) \rho^{\nu-\alpha} + \cdots$,  $a_\rho(\rho^\alpha u, \rho) = (\nu B(u) - \alpha u B'(u)) \rho^{\nu-1} + \cdots$, $\hat{\Theta}_x(\rho^\alpha u, \rho) = P'(u) \rho^{r-1} + \cdots$, and $V_\rho(\rho^\alpha u, \rho) = ((r-1) Q(u)- \alpha u Q'(u)) \rho^{r-2} + \cdots$, where the dots denote higher order terms in $\rho$. In summary, we get $N(\rho^\alpha u, \rho) = \mathcal{N}(u) \rho^{r +\nu -1} + \cdots$ with $\mathcal{N} = P B' - B P' +(m+\nu-r) B Q - \alpha u (B Q'-Q B')$. Thus, since $T(\rho^\alpha u, \rho) = D(u) \rho^{r+\alpha-1}+ \cdots$, we find
$$
\lim_{\rho \to 0} \frac{\rho^\alpha N(\rho^\alpha u, \rho)}{a(\rho^\alpha u, \rho) \, T(\rho^\alpha u, \rho)} = \frac{\mathcal{N}(u)}{B(u) D(u)} = \Xi(u),
$$
where $\Xi(u)$ is defined in \eqref{eq:Xi1}. Since for each fixed $u$, $\rho(u; \rho_0) \to 0$ when $\rho_0 \to 0^+$ we have proved that  the integrand in \eqref{int-gamma} converges pointwise to $\Xi(u)$. We conclude, by dominated convergence and \eqref{eq:majorant}, that
\[
 \lim_{\rho_0 \to 0^+} \int_{\gamma_{\rho_0}} E\, d t = \int_{\mathbb{R}} \Xi(u)\, d u,
\]
and the limiting integral is absolutely convergent. Substitution into \eqref{ec-LogH-C} gives \eqref{eq:Cconnection} and this completes the proof.
\end{proof}

\subsection{Proof of Proposition \ref{prop-perfil}}

\begin{proof}
We assume all the hypothesis of Theorem \ref{teo-Fund-G-jet}. In particular, from \eqref{eq:positiveDB}, $B(u) > 0$ and $D(u) > 0$ for every $u\in\mathbb R$.

Assume that a function $h(u)$ exists such that $\mathcal{G}(\varphi^* + \rho^\alpha u, \rho) = h(u) \rho^{m + \alpha \kappa} + o(\rho^{m + \alpha \kappa})$. We are going to see that it is a solution of the first order differential equation \eqref{profile-edo}. Indeed, after the blow-up $x = \rho^\alpha u$, the operators $\partial_\rho$ and $\partial_\varphi$ are transformed into $(\partial u/ \partial \rho) \partial_u +  \partial_\rho = -\alpha  u/\rho \partial_u + \partial_\rho$ and $(\partial u/ \partial \theta) \partial_u  = 1/\rho^\alpha \partial_u$, respectively. Then the function $\bar{G}(u, \rho) := \mathcal{G}(\varphi^* + \rho^\alpha u, \rho)$ satisfies the partial differential equation
$$
\bar{\Theta}(u, \rho) \frac{1}{\rho^\alpha} \partial_u \bar{G} + \bar{R}(u, \rho) \left(-\alpha  \frac{u}{\rho} \partial_u + \partial_\rho \right) \bar{G} = \overline{{\rm div} \mathcal{Z}}(u, \rho)  \, \bar{G},
$$
where $\bar{\Theta}(u, \rho) := \Theta(\varphi^* +\rho^\alpha u, \rho) = \rho^{r+\alpha-1} P(u) + o(\rho^{r+\alpha-1})$, $\bar{R}(u, \rho) := R(\varphi^* +\rho^\alpha u, \rho) = \rho^r Q(u) + o(\rho^{r})$ and $\overline{{\rm div} \mathcal{Z}}(u, \rho) := {\rm div} \mathcal{Z}(\varphi^* +\rho^\alpha u, \rho) =(P'(u)+ r Q(u)-\alpha u Q'(u)) \rho^{r-1} + o(\rho^{r-1})$. Inserting $\bar{G}(u, \rho) = h(u) \rho^{m + \alpha \kappa} + o(\rho^{m + \alpha \kappa})$ into the former equation and equating the coefficients of the minimum power of $\rho$ (that is $\rho^{m+\alpha \kappa+r-1}$) we get equation \eqref{profile-edo}.

We recall that $h(u) \neq 0$ for every $u \in \mathbb{R}$ since it is a nontrivial solution of an homogeneous equation. A direct algebraic manipulation on equation \eqref{profile-edo} gives that $h'/h = E - \Xi$, where the function $\Xi(u)$ is defined in \eqref{eq:Xi1} and we define
$$
E(u) :=\frac{B'}{B}+(\nu-\alpha\kappa)\frac{Q}{D}.
$$
Integrating the former differential equation over the interval $[-R, R]$ gives
\[
\log\left|\frac{h(R)}{h(-R)}\right| = \int_{-R}^{R}E(u)\,du -\int_{-R}^{R}\Xi(u)\,du.
\]
As $R\to\infty$, the left-hand side tends to $\log|L|$. The first integral in the right-hand side tends to zero by the compatibility condition hypothesis and the second integral converges to $\log(C_+/C_-)$ by \eqref{eq:Cconnection}. This proves the result.
\end{proof}

Sufficient conditions for the existence of $L$ are given below.

\begin{proposition} \label{prop-perfil-2}
Let $h(u)$ be a solution of \eqref{profile-edo} that we write as $D(u) h'(u) = \mathcal{A}(u) h(u)$. If the limits
\[
L_+=\lim_{u\to+\infty}\frac{h(u)}{|u|^\kappa},
\qquad
L_-=\lim_{u\to-\infty}\frac{h(u)}{|u|^\kappa},
\]
exist and are non-zero then $L = L_+/L_-$. Moreover, a sufficient condition for the existence of $L_{\pm} \neq 0$ is that $\mathcal{A}(u)/D(u)$ has at most a simple pole at infinity, that is, it has the asymptotic behaviour
\begin{equation}\label{integ-suf-cond}
\frac{\mathcal{A}(u)}{D(u)} = \frac{\kappa}{u}+O(|u|^{-1-\varepsilon})\quad(u\to\pm\infty)
\end{equation}
for some $\varepsilon>0$.
\end{proposition}
\begin{proof}
We assume $L_{\pm}$ exist and are non-zero. Then
\[
L = \lim_{u \to \infty} \frac{h(u)}{h(-u)} = \lim_{u \to \infty} \frac{h(u)/u^\kappa}{h(-u)/u^\kappa} = \frac{L_+}{L_-}.
\]

Now we shall prove the second part. Since $h'/h = \mathcal{A}/D$, we obtain
$$
\frac{d}{d u} \log\left| \frac{h(u)}{|u|^\kappa} \right| = \frac{\mathcal{A}(u)}{D(u)} - \frac{\kappa}{u}.
$$
If we assume integrable this right-hand side in $[R, \infty)$ and in $(-\infty, R]$ for any $R>0$ it follows that the function $h(u)/|u|^\kappa$ must have nonzero limit when $u\to\pm\infty$ and this proves the existence of $L_{\pm} \neq 0$. The proof finishes noticing that \eqref{integ-suf-cond} is a sufficient condition for the former integrability assumption.
\end{proof}

\begin{remark}
{\rm In $\Lambda \setminus \Lambda_{pq}$ with $0 \not\in \Omega_{pq}$ there are initial conditions $\mathcal{G}(0, \rho)$ analytic at $\rho=0$ yielding $\left. \mathcal{G}(\varphi, \rho) \right|_{I_{pq}} = g_m(\varphi) \, \rho^m + \cdots$ and $g_m(\varphi) \not\equiv 0$ with both $m=1$ and also $m\geq 2$. Actually, taking into account that $\hat{g}_1(\varphi; c_1) = c_1 G_0(\varphi)$ with $G_0(0) \neq 0$ it follows that either $m=1$ or $m \geq 2$ according to whether $c_1 \neq 0$ or $c_1 = 0$ that corresponds to whether $\mathcal{G}(0, \rho) = c_1 \rho + \cdots$ with $c_1 \neq 0$ or $\mathcal{G}(0, \rho) = O(\rho^2)$.}
\end{remark}

\subsection{Proof of Theorem \ref{cor-calculo-eta1}}

\begin{proof}
We take $\left. \mathcal{G}(\varphi, \rho) \right|_{I_{pq}} = g_m(\varphi) \, \rho^m + \cdots$ with $m \geq 2$ and $g_m(\varphi) \not\equiv 0$. Then equation
\begin{equation}\label{eta1-general}
g_m(2 \pi) \eta_1^{m-1} - g_m(0) = 0
\end{equation}
is just the relation obtained after introducing this $\left. \mathcal{G}(\varphi, \rho) \right|_{I_{pq}}$ into \eqref{G-fundamental2} and equating coefficients of the monomial $\rho^m$.

Let $0 \not\in \Omega_{pq} = \{ \varphi_1^*, \ldots, \varphi_s^* \}$ with the angles sorted in increasing order. We define the connected components $I_j$ of $I_{pq}$ such that $I_{pq} = \cup_j I_j$ with $I_1 = [0, \varphi_1^*)$, $I_j = (\varphi_{j-1}^*, \varphi_{j}^*)$ for $j=2, \ldots, s-1$ and $I_s = (\varphi_s^*, 2\pi]$ where $s = \# \Omega_{pq}$. Then we know that $g_m(\varphi) = C_{-}^{[j]} \hat{g}_m(\varphi)$ when $\varphi \in I_j$ while $g_m(\varphi) = C_{+}^{[j]} \hat{g}_m(\varphi)$ when $\varphi \in I_{j+1}$ where, of course,  $C_{+}^{[j]} = C_{-}^{[j+1]}$ for consistency.

Since $\eta_1 > 0$, the solution of \eqref{eta1-general} is
\begin{equation}\label{eta1-provisional}
\eta_1 = \left( \frac{g_m(0)}{g_m(2 \pi)} \right)^{\frac{1}{m-1}} = \left( \frac{C_{-}^{[1]} \hat{g}_m(0)}{C_{+}^{[s]} \hat{g}_m(2 \pi)} \right)^{\frac{1}{m-1}}.
\end{equation}
Since $\hat{g}_m$ satisfies \eqref{EDO-GM} and it is analytic at $\Omega_{pq}$ we have
$$
\hat{g}_m(\varphi) = c \, \exp\left( - PV \int_0^\varphi F(s) \, ds \right),
$$
for some constant $c \neq 0$ and $F := ((m-1) R_1 - G'_0)/G_0$. This Cauchy principal value exists because the analyticity of $\hat{g}_m$ at $\Omega_{pq}$ is equivalent to $F$ have at most a simple pole at $\Omega_{pq}$. Therefore
$$
\left( \frac{\hat{g}_m(0)}{\hat{g}_m(2 \pi)} \right)^{\frac{1}{m-1}} = \exp\left( \frac{1}{m-1} PV \int_0^{2 \pi} F(\theta) \, d \theta \right).
$$
We notice that $PV \int_0^{2 \pi} G'_0/G_0 \, d \varphi = 0$ since $G_0$ is $2 \pi$-periodic, $C^1([0, 2 \pi])$, and has a finite number of zeros in $[0, 2 \pi]$. Therefore, from the definition of $F$, one has
$$
PV \int_0^{2 \pi} F(\theta) \, d \theta = (m-1) \, PV \int_0^{2 \pi} \frac{R_{1}(\varphi)}{G_{0}(\varphi)} d \varphi.
$$
Then the expression of $\eta_1$ given in \eqref{eta1-provisional} becomes
$$
\eta_1 = \left( \frac{C_{-}^{[1]}}{C_{+}^{[s]}} \right)^{\frac{1}{m-1}} \, \exp\left( PV \int_0^{2 \pi} \frac{R_{1}(\varphi)}{G_{0}(\varphi)} d \varphi \right).
$$
We define
$$
C_m := \frac{C_{-}^{[1]}}{C_{+}^{[s]}} = \prod_{j=1}^s  \frac{C_{-}^{[j]}}{C_{+}^{[j]}}  = \prod_{j=1}^s \exp\left( \int_{\mathbb{R}} \Xi_j(u) \, du \right),
$$
where in the last step we have used Theorem \ref{teo-Fund-G-jet} applied to each $\varphi_j^* \in \Omega_{pq}$. This finishes the proof.
\end{proof}

\subsection{Proof of Lemma \ref{lemma-edo-homogenea-formal}}

\begin{proof}
First we place the singularity $x^*$ at the origin by the translation $x \mapsto x-x^*$ that does not affect the homogeneous linear structure of the differential equation. We still denote, for simplicity, $A$ and $B$ the resulting coefficients.

Let $m$ be the multiplicity of $A$ at $x=0$, that is, $A(x) = x^m u_0(x)$ and we assume there is a formal solution $y(x) = x^r u_1(x)$ with $r \in \mathbb{N}$. Here $u_j$ is an invertible formal series centeres at the origin, hence $u_j(0) \neq 0$ for $j=0,1$. Inserting $A$ and $y$ in the differential equation $A(x) y' + B(x) y = 0$ and dividing by $x^r$ we get $r x^{m-1} u_0 u_1 + x^m u_0 u_1' + B u_1 = 0$. Since $u$ is invertible we solve for $B =- r x^{m-1} u_0 - x^m u_0 u_1'/u_1$. Since $A' = m x^{m-1} u_0 + x^m u_0'$, we write $B + \frac{r}{m} A' = x^m \left(\frac{r}{m} u_0' - u_0 \frac{u_1'}{u_1} \right)$. Defining the analytic function $F(x) = \frac{r}{m} \frac{u_0'}{u_0} - \frac{u_1'}{u_1}$ we write
\begin{equation}\label{def-anal-A}
B = - \frac{r}{m} A' + A F.
\end{equation}
Since $A'/A = m/x + \cdots$, equation \eqref{def-anal-A} implies that $-B/A = r/x + \cdots$ where the dots means analytic higher order terms. So $B/A$ has a simple pole at the origin and $r = {\rm Res}(-B/A, 0)$.

All these conditions are necessary for the existence of the nontrivial formal solution $y(x)$. To prove its convergence we observe that the function $u_1 = y(x)/x^r$ is a solution of the analytic linear differential equation $u_1' + F u = 0$ that is non-singular at $x=0$ and consequently, given an initial condition $u_1(0) \neq 0$, produces a unique function $u_1$ analytic at $x=0$, hence $y$ is also analytic at $x=0$. This proves the first part.

To prove the second part, we assume now that $-B/A =  r/x + h(x)$ with $r$ a nonnegative integer and $h$ an analytic function at $x=0$. Then we get a primitive
$$
\mathcal{P}(x) = \int -B(x)/A(x) \, dx = \log(|x|^r) + H(x)
$$
with $H$ an analytic function at $x=0$. Notice that $\mathcal{P}$ is not continuous at $x=0$. Therefore, for any $x_0<0$ and $x>0$, any non-special solution has the form
\begin{eqnarray*}
y(x) &=& y(x_0) \exp\left( PV \int_{x_0}^x -B(s)/A(s) \, ds\right)  \\
 &=&  y(x_0) \exp\left( \lim_{\varepsilon \to 0^+} (\mathcal{P}(-\varepsilon) - \mathcal{P}(x_0) + \mathcal{P}(x)-\mathcal{P}(\varepsilon) )  \right) = |x|^r \hat{H}(x)
\end{eqnarray*}
with $\hat{H}$ an analytic function at $x=0$. We conclude that $y(x)$ will be an analytic solution at $x=0$ provided that $r \in 2 \mathbb{N}_0$.
\end{proof}

\subsection{Proof of Lemma \ref{lemma-multipl}}

\begin{proof}
Given weights $(p,q) \in W(\mathbf{N}(\mathcal{X}))$, we consider the $(p, q)$-quasihomogeneous expansion $\mathcal{X} = \mathcal{X}_r + \cdots$ where the leading $(p, q)$-quasihomogeneous vector field $\mathcal{X}_{r} = P_{p+r}(x,y) \partial_x + Q_{q+r}(x,y) \partial_y$. Therefore the polynomials $P_{p+r}$ and $Q_{q+r}$ are $(p, q)$-quasihomogeneous and, consequently, they have a factorization in $\mathbb{R}[x,y]$ of the form $\prod_j f_j^{m_j} \, \prod_j g_j^{n_j}$ with $f_j(x,y)$ either $x$ or $y$ or $y^p- \lambda_j x^q$ with $\lambda_j \neq 0$ while $g_j(x,y) = (y^p- a_j x^q)^2 + b_j^2 x^{2 q}$ with $b_j \neq 0$.

In $(p,q)$-weighted polar coordinates $(x,y) \mapsto (\varphi, \rho)$ given by
\begin{equation}\label{wpbu2}
x= \rho^p \cos\varphi, \ \ y = \rho^q \sin\varphi,
\end{equation}
the vector field $\mathcal{X}_r$ is transformed into $\mathcal{Z}_r = G_0(\varphi) \partial_\varphi + R_1(\varphi) \, \rho \, \partial_\rho$
after removing the common factor $\rho^{r} / D(\varphi) > 0$ with $D(\varphi) = p \cos^2\varphi + q \sin^2\varphi > 0$. Here
\begin{eqnarray}
R_1(\varphi) &=& P_{p+r}(\cos\varphi, \sin\varphi) \cos\varphi + Q_{q+r}(\cos\varphi, \sin\varphi) \sin\varphi, \label{def-R1}\\
G_0(\varphi) &=& p \, Q_{q+r}(\cos\varphi, \sin\varphi) \cos\varphi -q \, P_{p+r}(\cos\varphi, \sin\varphi) \sin\varphi \label{def-G01}.
\end{eqnarray}
From \eqref{def-G01} we see that the common factors of type $f_j(x,y)$ in $P_{p+r}(x,y)$ and $Q_{q+r}(x,y)$ give rise to factors $f_j(\cos\varphi, \sin\varphi)$ in $G_0$ such that $f_j(\cos\varphi_j^*, \sin\varphi_j^*) = 0$ for each $\varphi_j^* \in \Omega_{pq}$. Notice that $g_j(\cos\varphi, \sin\varphi) \neq 0$ for any $\varphi \in \mathbb{S}^1$.

In summary, generically if no cancelations happen, the factors  $\sin^{2k}\varphi$, $\cos^{2k}\varphi$ and $(\sin^p\varphi-\lambda \cos^q\varphi)^{2k}$ of $G_0$ come from the common factors $\sin^{2k-1}\varphi$, $\cos^{2k-1}\varphi$ and $(\sin^p\varphi-\lambda \cos^q\varphi)^{2k}$ of the trigonometric polynomials $P_{p+r}(\cos\varphi, \sin\varphi)$ and $Q_{q+r}(\cos\varphi, \sin\varphi)$, respectively. They turn via \eqref{def-R1} also in factors of $R_1$ proving the lemma.
\end{proof}

%%%%%%%%%%%%%%%%%%%%%%%%%%%%%%%%%%%%%%%%%%%%%%%%%%%%%

\end{document}